\documentclass[12 pt]{article}
\usepackage{amssymb,amsmath,amsfonts,amsthm}
\usepackage{xcolor}
\theoremstyle{plain}

\theoremstyle{remark}

\newcommand\BBP{{\mathbb {P}}}

\newcommand\BBE{{\mathbb {E}}}

\newcommand\E{{\mathbb {E}}}

\newtheorem {Theorem}{Theorem}[section]
\newtheorem {Proposition}{Proposition}[section]
\newtheorem {Corollary}{Corollary}[section]

\theoremstyle{definition}

\newtheorem {Fact}{Fact}[section]

\newcommand\beq{\begin{equation}}
\newcommand\eeq{\end{equation}}

\def\BBE{\mathbb{E}}
\def\BBP{\mathbb{P}}

\begin{document}

\title{Rates in the central limit theorem for random projections of Martingales}

\author{J. Dedecker\footnote{J\'er\^ome Dedecker, Universit\'e de Paris, CNRS, MAP5, UMR 8145,
45 rue des  Saints-P\`eres,
F-75006 Paris, France.}, F. Merlev\`ede \footnote{Florence Merlev\`ede, LAMA,  Univ Gustave Eiffel, Univ Paris Est Cr\'eteil, UMR 8050 CNRS,  \  F-77454 Marne-La-Vall\'ee, France.}
and 
M. Peligrad \footnote{Magda Peligrad, Department of Mathematical Sciences, University of Cincinnati, PO Box 210025, Cincinnati, Oh 45221-0025, USA.}}

\maketitle

\begin{abstract}
In this paper, we consider partial sums of martingale differences weighted
by random variables drawn uniformly on the sphere, and globally independent of the martingale differences. Combining Lindeberg's
method and a series of arguments due to Bobkov, Chistyakov and G\"{o}tze, we
show that the Kolmogorov distance between the distribution of these weighted
sums and the limiting Gaussian is ``super-fast'' of order $(\log n)^{2}/n$,
under conditions allowing us to control the higher-order conditional moments
of the martingale differences. We also show that the same rate is achieved if we consider a  quantity very close to  these weighted
sums, and  give an application of this result to
the least squares estimator of the slope in the linear model with Gaussian
design.

\end{abstract}



\section{Introduction}

\setcounter{equation}{0}

Let $(X_{k})_{k\in {\mathbb{Z}}}$ be a strictly stationary process on a
probability space $(\Omega ,\mathcal{F},{\mathbb{P}})$ such that $\BBE(X_{0})=0$ and $%
\BBE(X_{0}^{2})=1$. Consider the weighted sums 
\begin{equation*}
\langle X,\theta \rangle :=\theta _{1}X_{1}+\cdots +\theta_{n}X_{n}=: S_{n}(\theta )\,,
\end{equation*}%
where $X =(X _{1},\dots ,X _{n})^t$  and $\theta =(\theta _{1},\dots ,\theta _{n})^t$ is defined on $(\Omega ,\mathcal{F},{\mathbb{P}})$, is independent of $%
(X_{k})_{k\in {\mathbb{Z}}}$ and has uniform distribution $\mu _{n-1}$ on
the unit sphere $S^{n-1}$ of ${\mathbb{R}}^{n}$ ($n\geq 2$). Let $({\mathcal{%
F}}_{i})_{i\in {\mathbb{Z}}}$ be a non decreasing stationary filtration in $%
\mathcal{F}$ such that $X_{0}$ is ${\mathcal{F}}_{0}$-adapted. In this paper, we shall often use the notation $\BBE_{\ell}(X)$ to mean 
$\BBE(X| {\mathcal{F}}_{\ell})$. Assume that $\BBE_{i-1}(X_{i})=0$ ${\mathbb{P}}$-almost surely: in other words $(X_{k})_{k\in 
{\mathbb{Z}}}$ is a sequence of martingale differences adapted to $({%
\mathcal{F}}_{k})_{k \in {\mathbb{Z}}}$. By the properties of the uniform distribution on the sphere, we have that  
\begin{equation*}
\max_{1\leq i\leq n}|\theta _{i}|\rightarrow 0  , \ \ \text{ ${\mathbb{P}}$-a.s., as $n \rightarrow \infty$.}
\end{equation*} 
Therefore, if we assume in addition that $\mathbb{E}(X_{0}^{2}|{\mathcal{F}}_{-\infty})=1$ almost surely, according to Hannan \cite[p. 284]{Ha79}, $S_{n}(\theta )$ converges in distribution to a standard
Gaussian random variable. Actually \cite[Theorem 1]{Ha79} implies the following
conditional version of the central limit theorem (CLT): 
\[
\kappa_{\theta} \big ( P_{S_n(\theta)}, P_G)  \rightarrow 0 \ \ \text{ ${\mathbb{P}}$-a.s., as $n \rightarrow \infty$,}
\]
where $G \sim {\mathcal N}(0,1)$ and
\[
\kappa_{\theta} \big ( P_{S_n(\theta)}, P_G)  = \sup_{t \in {\mathbb R}} \big | \BBP_{| \theta} ( S_n(\theta) \leq t) -  \BBP (G\leq t)  \big |  \, .
\] 
Above, the notation $ \BBP_{| \theta}$ (resp. $ \BBE_{| \theta}$) means  the conditional probability (resp. the conditional expectation) with respect to  $\theta$.

In this paper, we are interested in rates in this CLT in terms of the
Kolmogorov distance.

When $(X_{k})_{k\in {\mathbb{Z}}}$ forms a sequence of independent centered
random variables in ${\mathbb{L}}^{4}$ with variance one, from Corollary 3.2
in Klartag and Sodin \cite{KS}, it follows that 
\begin{equation}
{\mathbb{E}}\Big (\kappa _{\theta }(P_{S_{n}(\theta )},P_{G})\Big )\leq 
\frac{cM_{4,n}}{n}\,\text{ where }\,M_{4,n}=\frac{1}{n}\sum_{k=1}^{n}{%
\mathbb{E}}(X_{k}^{4})\,.  \label{resultindependent}
\end{equation}%
This proves that, when $M_{4,n}$ is uniformly bounded (which is the case in the setting of
independent and identically distributed (iid) r.v.'s), projecting the
variables on the sphere allows to derive a much faster rate than in the
usual Berry-Esseen theorem, where the rate is $1/\sqrt{n}$. In a recent
paper, Bobkov et al. \cite{BCG20b} have extended this interesting phenomenon
to isotropic random vectors (meaning that the coordinates are uncorrelated  with variance one) having a symmetric distribution, and under a suitable second order
correlation condition. They obtained a similar $1/n$-rate modulo a
logarithmic factor. More precisely, their second order correlation condition
reads as: there exists a constant $\Lambda $ such that, for any $n \geq 1$ and   any collection $%
a_{ij}\in {\mathbb{R}}$, 
\begin{equation}
\mathrm{Var}\Big (\sum_{i,j=1}^{n}a_{ij}X_{i}X_{j}\Big )\leq \Lambda
\sum_{i,j=1}^{n}a_{ij}^{2}\,.  \label{condsecondorder}
\end{equation}%
Theorem 1.1 in \cite{BCG20b}  asserts that if $(X_{k})_{k\in {\mathbb{Z}}%
}$ is a sequence of uncorrelated centered random variables with variance
one, satisfying \eqref{condsecondorder} and such that $(X_{1},\cdots ,X_{n})$
has a symmetric distribution then 
\begin{equation}
{\mathbb{E}}\Big (\kappa _{\theta }(P_{S_{n}(\theta )},P_{G})\Big )\leq 
\frac{c\log n}{n}\Lambda \,.  \label{resultsecondorder}
\end{equation}%
As shown in \cite{BCG20b}, condition \eqref{condsecondorder} can be verified
for random vectors which satisfy a Poincar\'{e}-type inequality with
positive constant (see \cite[Proposition 3.4]{BCG20b}). Moreover, if we
assume that $(X_{k})_{k\in {\mathbb{Z}}}$ is a sequence of martingale
differences such that $\sup_{i\geq 1}{\mathbb{E}}(X_{i}^{4})<\infty $, one
can check that condition \eqref{condsecondorder} is satisfied provided that 
\begin{equation}
\sum_{k\geq 1} k {\tilde{\gamma}}(k)<\infty \,,  \label{condsymmart}
\end{equation}%
with ${\tilde{\gamma}} (k) =\max (\gamma _{2,2}(k),\gamma _{1,3}(k))$ where 
\begin{equation}
\gamma _{2,2}(k)=\sup_{\ell \geq u\geq 0}\Vert X_{u}X_{\ell }({\mathbb{E}}%
_{\ell }(X_{k+\ell }^{2})-{\mathbb{E}}(X_{k+\ell }^{2}))\Vert _{1}
\label{defgamma22mart}
\end{equation}%
and 
\begin{equation}
\gamma _{1,3}(k)=\sup_{\ell ,v\geq 0}\Vert X_{\ell }\big ({\mathbb{E}}_{\ell
}(X_{k+\ell }X_{k+v+\ell }^{2})-{\mathbb{E}}(X_{k+\ell }X_{k+v+\ell }^{2})%
\big )\Vert _{1}\,.  \label{defgamma13mart}
\end{equation}%
For instance, if $(X_{k})_{k\in {\mathbb{Z}}}$ is additionally strictly
stationary and strongly mixing, condition \eqref{condsymmart} is satisfied
provided that $\sum_{k\geq 1} k \int_{0}^{\alpha _{2}(k)}Q^{4}(u)du<\infty $
(see Section \ref{MC} for a definition of the coefficients $\alpha _{2}(k)$
and of the quantile function $Q$). The most stringent condition in the
assumptions made in \cite[Theorem 1.1]{BCG20b} is probably the fact that the
distribution of $(X_{1},\cdots ,X_{n})$ is assumed to be symmetric. In \cite[%
Chapter 17.4]{BCGbook} the authors consider the case of non-symmetric
distributions. Their Proposition 17.4.1 states that if $(X_{k})_{k\in {%
\mathbb{Z}}}$ is a sequence of uncorrelated centered random variables with
variance one, satisfying \eqref{condsecondorder} then we can still provide
an explicit bound for ${\mathbb{E}}\big (\kappa _{\theta } \big (%
P_{S_{n}(\theta )},P_{G}\big )\big )$ and  a certain term has to be added to the right-hand
side of \eqref{resultsecondorder}. This additional term is 
\begin{equation}
\Big (\frac{\log n}{n}\Big )^{1/4}\Big ({\mathbb{E}}\frac{\langle X,Y\rangle 
}{\sqrt{\Vert X\Vert_e^{2}+\Vert Y\Vert_e^{2}}}\Big )^{1/2}\,,
\label{additionalterm}
\end{equation}%
where $Y$ is an independent copy of $X$ and $\Vert X\Vert_e^{2} =\langle X, X \rangle $ denotes here and all along the paper, the euclidian norm of $X$. As proved in \cite[Chapter 17.5]%
{BCGbook}, the term \eqref{additionalterm} can be upper-bounded by 
$C(\log n)^{1/4}/n$ provided $(X_{1},\ldots ,X_{n})$ satisfies a Poincar\'{e}-type
inequality. For instance, when $n=1$, for this Poincar\'{e}-type inequality
to be satisfied it is necessary that ${\mathbb{P}}_{X_{1}}$ has an
absolutely continuous component. Now in case of random vectors ($n\geq 2$),
the required Poincar\'{e}-type inequality is quite complicated to obtain
except in the case where the random variables are independent with marginal
distributions satisfying a Poincar\'{e}-type inequality (see \cite[Pages
108-109]{BCGbook}).

The aim of this paper is to provide a new method  allowing us to show that sequences $%
(X_{k})_{k\in {\mathbb{Z}}}$ of martingale differences  satisfy an upper
bound of the type \eqref{resultindependent} (up to some logarithmic term) without requiring that the law
of the vector $(X_{1},\cdots ,X_{n})$ is symmetric, nor satisfies a Poincar\'{e}-type inequality. 
More precisely, as stated in Theorem \ref{BEmartthm} and proved in Section \ref%
{Section3}, an upper bound of the type \eqref{resultindependent} will be
achieved with the help of the Lindeberg method, where the random variables $%
\theta _{i}X_{i}$ will be replaced one by one by random variables $%
Y_{i}(\theta )$ taking only two values with some desired characteristics (in particular they are independent conditionally to $\theta$).
Then, following the approach of Klartag and Sodin  \cite{KS}, we shall compare the distribution of $\sum_{i=1}^{n}Y_{i}(\theta )$
with a normal distribution. We select this method because  Stein's method,
which is also successfully used in some situations to get sharp Berry-Esseen
bounds, would require  strong assumptions on the conditional moments of the martingale differences (see for instance \cite{Rol18}). 

The rate we achieve is ``super fast"
of order $(\log n)^{2}/n.$ It should be mentioned that, in the absence of
the randomization considered in this section and in the presence of
dependence, it is very rare to achieve a Berry-Esseen upper bound of order $%
n^{-1/2}$ for strictly stationary sequences. In particular, El Machouri and Voln\'y \cite{EMV}  exhibited an example showing that the rate of convergence in the central limit theorem in terms of the Kolmogorov distance can be arbitrarily slow for strictly stationary, strong mixing and bounded martingale difference sequences (note that for this example the coefficients $\gamma _{2,2}(k) $ and $\gamma _{1,3}(k)$ 
defined in \eqref{defgamma22mart} and \eqref{defgamma13mart} converge to zero as $k \rightarrow \infty$).  To our knowledge, the best known Berry-Esseen
bound for strictly stationary and bounded martingale difference sequences is $O (n^{-1/3})$ under the condition $\sum_{k >0} k \theta_{X,3,4}(k) < \infty$ where $ \theta_{X,3,4}(k)$ is a slightly more restrictive coefficient than ${\tilde{\gamma}} (k) =\max (\gamma _{2,2}(k),\gamma _{1,3}(k))$ (see Theorem 2.7 in \cite{DMR23} and its implication on the Berry-Esseen type estimates as described in \cite[Corollary 2.4]{DMR23}).  To be complete, in the non stationary setting, even in case of  constant conditional variance and 
moments of order $3$ uniformly bounded, the upper bound in the classical Berry-Esseen inequality cannot be better
than $n^{-1/4}$ without additional assumptions (see Example 1 in  \cite{Botthausen}).

In a Berry-Esseen bound the constants are important. However, in this paper,
due to the difficulty and complexity of the problem, we shall not compute the
constants exactly. For any two positive sequences of random variables $a_{n}$
and $b_{n}$ we shall often use the notation  $a_{n}\ll b_{n}$  to mean that there is a 
positive finite constant $c$  such that $a_{n}\leq cb_{n}$ for all $n$. The constant $c$ is allowed to depend on the moments of $X_0$ and on some quantities involving coefficients such as 
$({\tilde{\gamma}} (k) )_{k \in {\mathbb Z}}$ but not on $n$. Also, all along the paper, the vectors will be in column form, ${\bf 1}_n =(1, 1, \cdots, 1)^t$ will be the unit vector of size $n$,  $I_n$ will designate the identity matrix of order $n$ and $J_n= {\bf 1}_n{\bf 1}_n^t  $ will be the all-ones square matrix  of order $n$.

The paper is organized as follows. In Section \ref{Section2.1} we present our main 
result concerning the rate in the CLT for ${ S}_n(\theta)$ in the case where $(X_k)_{k \in {\mathbb Z}}$ is a strictly stationary sequence of martingale differences  (Theorem \ref{BEmartthm}). We shall also prove that the same rate is achieved when we consider a quantity very close to $S_n(\theta)$, namely the quantity ${\tilde S}_n(\theta)= \langle X, A\theta \Vert A \theta \Vert_e^{-1} \rangle $ where $A=I_n - n^{-1}J_n$.  Applications are given in Section \ref{SectionAppli}. In particular, Sections \ref{MC} and \ref{ARCH} are devoted to applications to martingale differences satisfying
a mixing-type condition, harmonic functions of Markov chains and ARCH($\infty $)
models, whereas,  in Section \ref{Section4}, we apply our result to the ordinary least square estimator of the slope in
the linear regression model with Gaussian design and martingale differences
errors. The proof of our main result result is given in Section \ref{Section3}. 

\section{Normal approximation for weighted sums of martingale differences} \label{Section2}

\setcounter{equation}{0}

\subsection{Main Result} \label{Section2.1}

In what follows we assume that $(X_{i})_{i\in {\mathbb{Z}}}$ is a strictly stationary
sequence of martingale differences. Assume moreover that $\Vert X_{0}\Vert
_{4}<\infty $ and $\Vert X_{0}\Vert _{2}=1$. Let us introduce the weak
dependence coefficients we will use in this paper. For any positive integer $%
v$, let 
\begin{equation*}
\gamma _{0,2}(v)=\Vert {\mathbb{E}}_{0}(X_{v}^{2})-{\mathbb{E}}%
(X_{v}^{2})\Vert _{1}\,,\,\gamma _{1,2}(v)=\Vert X_{0}({\mathbb{E}}%
_{0}(X_{v}^{2})-{\mathbb{E}}(X_{v}^{2}))\Vert _{1}\,.
\end{equation*}%
Recall also the coefficients $\gamma _{2,2}(v)$ and $\gamma _{1,3}(v)$ as
defined in the introduction that can be rewritten as follows in the strictly
stationary setting: 
\begin{align*}
\gamma _{2,2}(v)&=\sup_{\ell \geq 0}\Vert X_{0}X_{\ell }({\mathbb{E}}_{\ell
}(X_{v+\ell }^{2})-{\mathbb{E}}(X_{v+\ell }^{2}))\Vert _{1}\,, \\
\gamma
_{1,3}(v)&=\sup_{\ell \geq 0}\Vert X_{0}\big ({\mathbb{E}}_{0}(X_{v}X_{v+\ell
}^{2})-{\mathbb{E}}(X_{v}X_{v+\ell }^{2})\big )\Vert _{1}\,.
\end{align*}%
Define then 
\begin{equation}
\gamma (v) =  \max (\gamma _{0,2}(v),\gamma _{1,2}(v),\gamma
_{2,2}(v),\gamma _{1,3}(v))\,.  \label{defgamma}
\end{equation}%
Our general result for martingales is the following theorem:

\begin{Theorem}
\label{BEmartthm} Let $(X_{n})_{n\in \mathbb{Z}}$ be a strictly stationary sequence
of martingale differences in ${\mathbb{L}}^{4}$ such that ${\mathbb{E}}%
(X_{0}^{2})=1$. Let $(\gamma (k))_{k\geq 0}$ be the sequence of dependent
coefficients defined in \eqref{defgamma}. Assume that $\sum_{k\geq 1}k\gamma
(k)<\infty $. Then there exists a positive constant $C_1$ such that for any $n \geq 2$, 
\begin{equation}
{\mathbb{E}}\big (\kappa _{\theta }(P_{S_{n}(\theta )},P_{G})\big )\leq C_1 \frac{%
(\log n)^{2}}{n}\,.  \label{concluBEmartthm}
\end{equation}
Moreover, with the notation $ \displaystyle {\tilde S}_n ( \theta) = \langle X, A \theta   \Vert  A \theta \Vert^{-1}_e   \rangle$ where $A= I_n - n^{-1} J_n$, there exists a positive constant $C_2$ such that for any $n \geq 2$, 
\begin{equation}
{\mathbb{E}}\big (\kappa _{\theta }(P_{ {\tilde S}_{n}(\theta )},P_{G})\big )\leq C_2 \frac{%
(\log n)^{2}}{n}\,.  \label{concluBEmartthm2}
\end{equation}
\end{Theorem}

As we shall see in Subsection \ref{Section4}, the upper bound \eqref{concluBEmartthm2} is very useful in the context of linear regression with Gaussian designs since in this case the rate of convergence in the central limit theorem for the least square estimator of the slope  is reduced to the study of the asymptotic behavior of $ {\tilde S}_n ( \theta) $. Note that in the definition of $ {\tilde S}_n ( \theta) $, the  self normalized quantity $A \theta   \Vert  A \theta \Vert^{-1}_e$ has, by definition, its euclidian norm equals to one but is not anymore uniformly distributed on the sphere because 
$ \langle A \theta   \Vert  A \theta \Vert^{-1}_e  , {\mathbf 1}_n  \rangle =0$. Therefore the upper bound \eqref{concluBEmartthm2} is not a direct application of \eqref{concluBEmartthm} even if some of the arguments to prove both upper bounds are similar.  

Note also that the upper bounds \eqref{concluBEmartthm} and \eqref{concluBEmartthm2} hold if we replace the random vector $\theta=(\theta
_{1},\dots,\theta_{n})$ by ${\tilde \xi}= \Vert \xi \Vert_e^{-1} \xi$ where $\xi = (\xi_{1},\dots,\xi_{n})^t$ with $(\xi_i)_{1 \leq i \leq n}$  iid centered and standard Gaussian r.v.'s independent of $(X_{j})_{1\leq j\leq n}$. Indeed, it is well-known that ${\tilde \xi}$  has uniform distribution on the unit sphere $S^{n-1}$ of
${\mathbb{R}}^{n}$.

The proof of Theorem \ref{BEmartthm} is given in Section \ref{Section3}.

\subsection{Applications}

\label{SectionAppli}

\subsubsection{Martingale differences sequences and functions of Markov
chains}

\label{MC}

Recall that the strong mixing coefficient of Rosenblatt \cite{Ro} between
two $\sigma$-algebras ${\mathcal{A}}$ and ${\mathcal{B}}$ is defined by $%
\alpha ({\mathcal{A}}, {\mathcal{B}})=\sup \{ |{\mathbb{P}}(A \cap B)-{%
\mathbb{P}}(A){\mathbb{P}}(B)| \, : \, (A, B) \in {\mathcal{A}} \times {%
\mathcal{B}} \, \} $. For a strictly stationary sequence $(X_i)_{i\in {%
\mathbb{Z}}}$, let ${\mathcal{F}}_i =\sigma (X_k, k \leq i)$. Define the
mixing coefficients $\alpha_2(n)$ of the sequence $(X_i)_{i\in {\mathbb{Z}}} 
$ by 
\begin{equation*}
\alpha_2(n)= \sup_{\ell \geq 0}\alpha ({\mathcal{F}}_0, \sigma(X_n, X_{n +
\ell}))\,.
\end{equation*}
For the sake of brevity, let $Q=Q_{X_0}$ where $Q_{X_0} $ is the quantile
function of $X_0$, that is the generalized inverse of $t \mapsto {\mathbb{P}}
( |X_0| >t )$. The coefficients $(\gamma (n))_{n\geq 0}$ can be controlled with the help of the coefficients  $(\alpha_2 (n))_{n\geq 0}$  and the quantile function $Q$ by using some inequalities given in Rio's book (see \cite[Theorem 1.1 and Lemma 2.1]{Rio17}), as done for instance to get Inequality (6.75) in \cite{MPU19}.  Hence, applying Theorem \ref{BEmartthm}, the following result holds:

\begin{Corollary}
\label{corstrong} Let $(X_{n})_{n\in \mathbb{Z}}$ be a strictly stationary sequence
of martingale differences in ${\mathbb{L}}^{4}$ such that ${\mathbb{E}}%
(X_{0}^{2})=1$. Assume that 
\begin{equation}
\sum_{k\geq 1}k\int_{0}^{\alpha _{2}(k)}Q^{4}(u)du<\infty 
\label{concluBEmartmixing}
\end{equation}%
Then the conclusion of Theorem \ref{BEmartthm} holds.
\end{Corollary}

Examples of martingale differences that are additionally strongly mixing are
harmonic functions of a Harris recurrent Markov chain. Let us for instance
consider the following example as described in \cite[Section 4.1]{DMR09}:
Let $(Y_i)_{i \in {\mathbb{Z}}}$ be the homogeneous Markov chain with state
space ${\mathbb{Z}}$ described at page 320 in \cite{Da}. The transition
probabilities are  then given by $p_{n,n+1}=p_{-n,-n-1}=a_n$ for $n \geq 0$, $%
p_{n,0}=p_{-n,0}=1-a_n$ for $n > 0$, $p_{0,0}= 0$, $a_0 =1/2$ and $1/2 \leq
a_n < 1$ for $n \geq 1$. This chain is aperiodic and positively recurrent as soon
as $\sum_{n \geq 2}\Pi_{k=1}^{n-1}a_k < \infty$ and in that case the
stationary chain is strongly mixing in the sense of Rosenblatt \cite{Ro}. In addition, if we denote by 
$\tau =\inf\{n>0, X_n=0\}$, according to 
\cite[Theorem 2]{Bol80}, if ${\mathbb{E}} ( \tau^p | X_0=0) < \infty$ for some $p >2$, then 
$\sum_{k \geq 1} k^{p-2} \alpha_2(k) <
\infty$.

Denote by $K$ the Markov kernel of the chain $(Y_i)_{i \in {\mathbb{Z}}}$.
The functions $f$ such that $K(f)=0$ almost everywhere are obtained by
linear combinations of the two functions $f_1$ and $f_2$ given by $f_1(1)=1$%
, $f_1(-1)=-1$ and $f_1(n)= f_1(-n)= 0$ if $n \neq 1$, and $f_2(0)=1$, $%
f_2(1)=f_2(-1)=0$ and $f_2(n+1)=f_2(-n-1) = 1 - a_n^{-1}$ if $n >0$. Hence
the functions $f$ such that $K(f)=0$ are bounded.

From the above considerations, if $(X_i)_{i \in \mathbb{Z}}$ is defined by $X_i = f(Y_i)$ with $K(f)=0$,
then Corollary \ref{corstrong} applies if $\sum_{k \geq 1} k \alpha_2(k) <
\infty$, which in turn holds if  ${\mathbb{E}} ( \tau^3 | X_0=0) < \infty$.  But ${\mathbb P} (\tau = n | X_0=0)
=(1- a_n)\Pi_{i=1}^{n-1}a_i$ for $n\geq 2$. Consequently, if for some $%
\epsilon >0$, $\displaystyle a_i = 1 - \frac{1}{i} \Big ( 3 + \frac{1 + \epsilon}{\log i} %
\Big )$ for $i$ large enough, Corollary \ref{corstrong} applies.

\subsubsection{ARCH models}

\label{ARCH}

Theorem \ref{BEmartthm} applies to the case where $(X_{i})_{i\in {\mathbb{Z}}%
}$ has an ARCH($\infty $) structure as described by Giraitis \textit{et al.} 
\cite{Gi2000}, that is 
\begin{equation}
X_{n}=\sigma _{n}\eta _{n},\ \text{with}\ \sigma _{n}\in {\mathbb{R}}^{+}\ 
\text{that satisfies}\ \sigma _{n}^{2}=c+\sum_{j=1}^{\infty
}c_{j}X_{n-j}^{2}\,,  \label{defARCH}
\end{equation}%
where $(\eta _{n})_{n\in \mathbb{Z}}$ is a sequence of iid centered random
variables such that ${\mathbb{E}}(\eta _{0}^{2})=1$ and independent of ${\mathcal F}_{n-1}$, and where $c\geq 0$, $%
c_{j}\geq 0$, and $\sum_{j\geq 1}c_{j}<1$. Since ${\mathbb{E}}(\eta
_{0}^{2})=1$ and $\sum_{j\geq 1}c_{j}<1$, the unique stationary solution to (%
\ref{defARCH}) is given by Giraitis \textit{et al.} \cite{Gi2000}: 
\begin{equation}
\sigma _{n}^{2}=c+c\sum_{\ell =1}^{\infty }\sum_{j_{1},\dots ,j_{\ell
}=1}^{\infty }c_{j_{1}}\dots c_{j_{\ell }}\eta _{n-j_{1}}^{2}\dots \eta
_{n-(j_{1}+\dots +j_{\ell })}^{2}\,.  \label{solARCH}
\end{equation}%
Let $v^{2}={\mathbb{E}}(X_{0}^{2})$ and note that $v^{2}=c\big (1-\frac{{}}{%
{}}\sum_{j\geq 1}c_{j}\big )^{-1}$. Applying Theorem \ref{BEmartthm}, we get
the following result:

\begin{Corollary}
\label{corARCH} Let $p\in ]4,6]$. Assume that  $\Vert \eta _{0}\Vert
_{p}<\infty $,  $\sum_{j\geq 1}c_{j}<1$ and 
\begin{equation}
c_{j} \leq O( j^{-b})\ \text{for}\ b>1+2(p-2)/(p-4)   \,.  \label{condARCH}
\end{equation}%
Then, for any $n \geq 2$, 
\begin{equation*}
{\mathbb{E}}\big (\kappa _{\theta }(P_{S_{n}(\theta )},P_{G_{v^{2}}})\big )%
\ll \frac{(\log n)^{2}}{n}\,,
\end{equation*}%
where $G_{v^{2}}\sim {\mathcal{N}}(0,v^{2})$.
\end{Corollary}

\noindent \textbf{Proof of Corollary \ref{corARCH}.} Let $k \geq 2$ and $f$ be such that 
$
\sigma_k = f ( \eta_{k-1}, \ldots, \eta_{1} , \eta_0, \eta_{-1}, \ldots) $.  Let $(\eta^{\prime }_n)_{n \in \mathbb{Z}}$ be an independent copy of $%
(\eta_n)_{n \in \mathbb{Z}}$ and set 
\begin{equation*}
\sigma^*_k = f ( \eta_{k-1}, \ldots, \eta_{1} , \eta_0^{\prime },
\eta_{-1}^{\prime }, \ldots) \, .
\end{equation*}
Let also $X_k^*= \sigma^*_k \eta_k$, $k \geq 2$. Let us estimate the
coefficients $( \gamma (k) )_{k \geq 2}$. With this aim, we start by
noticing that, for any $k \geq 2$, 
\begin{equation*}
\gamma_{0,2} (k) \leq {\mathbb{E}} ( | X^2_k -X_k^{*2} | ) = {\mathbb{E}} (
\eta_k^2) {\mathbb{E}} ( | \sigma^2_k - \sigma_k^{*2} | ) = {\mathbb{E}} ( |
\sigma^2_k - \sigma_k^{*2} | ) =: \delta_k \, .
\end{equation*}
Now, according to \cite[Prop. 5.1]{CDT}, 
\begin{equation*}  
\delta_k \ll \inf_{1 \leq \ell \leq k} \big \{ \kappa^{k / \ell} + \sum_{i
\geq \ell +1} c_j \big \} \text{ where } \kappa = \sum_{i \geq 1} c_j  < 1\, ,
\end{equation*}
which implies, since we assumed that $c_{j}=O(j^{-b})$, 
\begin{equation}  \label{controldeltak}
\delta_k \ll   ( k^{-1} \log k)^{b-1} \, .
\end{equation}
Next, for any $k \geq 2$, 
\begin{align*}
\gamma_{1,2} (k) & \leq {\mathbb{E}} \big ( |X_0| | X^2_k -X_k^{*2} | \big ) %
= {\mathbb{E}} \big ( |X_0| | \sigma^2_k -\sigma_k^{*2} | \big ) \, , \\
\gamma_{2,2} (k) & \leq \sup_{\ell \geq 0}{\mathbb{E}} \big ( |X_0X_\ell| |
X^2_{k+ \ell} -X_{k+ \ell}^{*2} | \big ) = \sup_{\ell \geq 0}{\mathbb{E}} %
\big ( |X_0X_\ell| | \sigma^2_{k+ \ell} -\sigma_{k+ \ell}^{*2} | \big ) \, ,
\end{align*}
and, there exists a numerical constant $C$ such that, for any $k \geq 2$, 
\begin{equation*}
\gamma_{1,3} (k) \leq C \sup_{\ell \geq 0}{\mathbb{E}} \big ( \{ |X_\ell |^3
+ |X^*_\ell |^3 \} | X_{k} -X_{k}^{*} | \big ) \, .
\end{equation*}
Let us give an upper bound for the coefficients $\gamma_{1,3} (k) $, the
other coefficients being bounded by similar arguments.

For any $k \geq 2$, let $u(k):= \sup_{\ell}\Vert |X_{\ell}|^3 (X_k - X_k^*)
\Vert_1$. Let $M$ be a positive real. By stationarity and Cauchy-Schwarz's
inequality, note that 
\begin{multline*}
u(k) \leq M^{3-p/2} \sqrt{ {\mathbb{E}} ( |X_0|^p)} \sqrt{{\mathbb{E}} ( |
X_k - X_k^* |^2 ) } + M^{4-p} \sup_{\ell}\Vert |X_{\ell}|^{p-1} (X_k -
X_k^*) \Vert_1 \\
\leq M^{3-p/2} \sqrt{ {\mathbb{E}} ( |X_0|^p)} \sqrt{{\mathbb{E}} ( | X_k -
X_k^* |^2 ) } + 2 M^{4-p} {\mathbb{E}} ( |X_0|^p) \, .
\end{multline*}
Note that 
\begin{equation*}
{\mathbb{E}} ( | X_k - X_k^* |^2 ) = {\mathbb{E}} ( | \sigma_k - \sigma_k^*
|^2 ) \leq {\mathbb{E}} ( | \sigma^2_k - \sigma_k^{*2} | ) \, .
\end{equation*}
Therefore 
\begin{equation*}
u (k) \ll M^{3-p/2} \sqrt{\delta_k} + M^{4-p} \, .
\end{equation*}
The quantity $\sup_{\ell}\Vert |X^*_{\ell}|^3 (X_k - X_k^*) \Vert_1$ can be
bounded similarly. Selecting $M= \delta^{-1/(p-2)}_k$, it follows that 
\begin{equation*}
\gamma_{1,3} (k) \ll \delta^{(p-4)/(p-2)}_k \, .
\end{equation*}
Using similar arguments we infer that $\gamma_{1,2} (k) \ll
\delta^{(p-3)/(p-2)}_k$ and $\gamma_{2,2} (k) \ll \delta^{(p-4)/(p-2)}_k $.

So, overall, Theorem \ref{BEmartthm} applies if $\sum_{k \geq 1} k
\delta^{(p-4)/(p-2)}_k < \infty$ that clearly holds under \eqref{condARCH}
by taking into account \eqref{controldeltak}. \qed

\subsubsection{Application to linear regression with Gaussian design}

\label{Section4} 

Let us consider the following linear model 
\begin{equation*}
Y_i = \alpha + \beta Z_i + X_i \, , \, 1 \leq i \leq n \, ,
\end{equation*}
where $(X_i)_{i \in {\mathbb{Z}}}$ is a strictly stationary sequence of
martingale differences, and $(Z_i)_{i \in {\mathbb{Z}}}$ is a sequence of iid $ {%
\mathcal{N}} (\mu,\sigma^2)$-distributed random variables, which is independent of $(X_i)_{i \in {\mathbb{Z%
}}}$.

 As usual, the observations are $(Y_i,Z_i)_{1 \leq i \leq n}$ and the
aim is to estimate the unknown parameter $\beta$. The ordinary least squares (OLS)
estimator ${\hat \beta}$ of $\beta$ is then given by 
\begin{equation*}
{\hat \beta} = \sum_{i=1}^n \frac{ (Y_i - {\bar Y}_n) (Z_i - {\bar Z}_n)}{
\Vert Z -{\bar Z}_n {\mathbf 1}_n \Vert_e^2} \, ,
\end{equation*}
where the  notations ${\bar U}_n = n^{-1} \sum_{i=1}^n U_i$ , ${\mathbf 1}_n = (1,1, \cdots, 1)^t$ and $Z= (Z_1, \dots, Z_n)^t$ are used.  

Noting that 
$\Vert Z -{\bar Z}_n {\mathbf 1}_n  \Vert_e^2 = \sum_{i=1}^n (Z_i - {\bar Z}_n)^2$, the OLS estimator satisfies 
\begin{equation*}
{\hat \beta} - \beta =  \sum_{i=1}^n \frac{ (Z_i - {\bar Z}_n)}{ \Vert Z -{%
\bar Z}_n  {\mathbf 1}_n\Vert_e^2} X_i  \, .
\end{equation*}
Let 
\begin{equation*}
T_n:= \Vert Z -{\bar Z}_n   {\mathbf 1}_n \Vert_e \big ( {\hat \beta} - \beta \big ) = \sum_{i=1}^n \frac{ (Z_i - {\bar Z}_n)}{ \Vert Z -{%
\bar Z}_n   {\mathbf 1}_n \Vert_e} X_i \, .
\end{equation*}
Note that 
\begin{equation*}
T_n = \Vert {\tilde \xi} \Vert_e^{-1} \langle {\tilde \xi} , X
\rangle_{{\mathbb{R}}^n} 
\end{equation*}
where  ${%
\tilde \xi} = ( \xi_1 - {\bar \xi}_n, \cdots, \xi_n - {\bar \xi}_n)$ with $%
\xi_i = (Z_i- \mu)/ \sigma$.

\medskip

Since $\Vert { \xi} \Vert_e^{-1} {\xi}$  has the same law as $\theta$ which is uniformly distributed on the sphere ${\mathcal S}^{n-1}$, 
$\Vert {\tilde \xi} \Vert_e^{-1}  {\tilde \xi}$ has the same law as $ A\theta \Vert A \theta \Vert_e^{-1}$ where $A=I_n - n^{-1} J_n$. Therefore, as a direct consequence of the upper bound \eqref{concluBEmartthm2} 
of Theorem \ref{BEmartthm}, the following result holds: 

\begin{Corollary}
\label{BEmartthmregression} Let $(X_n)_{n \in \mathbb{Z}}$ be a strictly stationary sequence of martingale differences in ${\mathbb{L}}^4$ such that ${\mathbb{E}%
} (X_0^2) = 1$. Assume that $\sum_{v \geq 1} v \gamma (v) < \infty$ (where $%
\gamma (k)$ is defined in \eqref{defgamma}). Then, for any $n \geq 2$, 
\begin{equation*}
{\mathbb{E}} \big ( \sup_{t \in {\mathbb{R}}} \big | {\mathbb{P}}_{| Z} (T_n
\leq t) - {\mathbb{P}} (G\leq t) \big | \big ) \ll \frac{(\log n)^2 }{n} \, .
\end{equation*}
\end{Corollary}

\section{Proof of Theorem \ref{BEmartthm}}

\label{Section3}

\setcounter{equation}{0}

\subsection{Proof of the upper bound \eqref{concluBEmartthm}} Following \cite{BCG18,BCG20b}, we shall use as a starting point a  variant of the smoothing inequality, which is custom
built for the type of randomization we used. Let us give some details: the so-called Berry-Esseen smoothing inequality
(see e.g. \cite[Ineq. (3.13) p. 538]{Feller}) together with Lemma 5.2 in 
\cite{BCG18}, which takes advantage of the fact that the random variables
are projected on the sphere, imply that there exists a positive constant $c_1$ such that, for all $n \geq 1$ and all $T\geq T_{0} \geq 1$, 
\begin{multline} \label{firstinegalite}
c_1 {\mathbb{E}}\big (\kappa _{\theta }(P_{\langle X,\theta \rangle },P_{G})\big
)\leq \int_{0}^{T_{0}}{\mathbb{E}}\big (|f_{\theta }(t)-\mathrm{e}^{-t^{2}/2}|%
\big )\frac{dt}{t} \\ +\frac{\log (T/T_{0})}{n} \big ( \sigma^2_{4,n}  + m^2_{4,n}  \big ) +\mathrm{e}^{-T_{0}^{2}/16}+\frac{%
1}{T}\, ,
\end{multline}
where the following notations are used: $f_{\theta }(t)={\mathbb{E}}\big (\mathrm{e}^{it\langle X,\theta
\rangle }|\theta \big )$, 
\[
\sigma_{4,n} :=  \frac{1}{ \sqrt{n}}\Big \Vert \sum_{k=1}^{n} \big ((X_{k} )^2 - {\mathbb{E}} (X_{k} )^2 %
\big ) \Big \Vert_2 \  \text{ and } \  m_{4,n}   := \frac{1 }{ \sqrt{n} } \Big ( \E \Big | \sum_{k=1}^{n}  X_{k}  Y_{k}  \Big |^4 \Big )^{1/4} \, , 
\]
with $Y=(Y_{1} , \cdots, Y_{n} )$  an independent copy of $X$. Note first that if 
\begin{equation}
\sum_{k\geq 1} \big |  {\rm Cov}(X_{0}^{2},X_{k}^{2}) \big | <\infty ,  \label{cond-sigma}
\end{equation}%
then, there exists $c_2>0$ such that for all $n \geq 1$,  $\sigma^2_{4,n}   \leq c_2$. 
Next, by \cite[Corollary 2.3]{BCG18}, 
\[
m^2_{4,n}  \leq \sup_{\theta \in S^{n-1}} {\mathbb{E}}(|\langle X,\theta \rangle |^{4}) \, .
\]
Recall now the following version of the so-called Burkholder's inequality: Suppose that $(d_k)_{ k \geq 1 }$ is a sequence of martingales differences in ${\mathbb L}^p$ for $ p \in [2, \infty[$.  Then, for any positive integer $n$, 
\begin{equation} \label{BI}
\Big \Vert \sum_{i=1}^n d_i  \Big \Vert^2_p  \leq (p-1)     \sum_{i=1}^n  \Vert d_i  \Vert^{2}_{p} \, .
\end{equation}
The constant $p-1$ is derived in \cite{Rio09}.  Applying Inequality \eqref{BI} with $p =4$,  we derive that,  for  any fixed point $\theta $ on the sphere $S^{n-1}$, 
\[
 \Big (  {\mathbb{E}}(|\langle X,\theta \rangle |^{4})  \Big )^{1/2} \leq 3  \sum_{i=1}^n  \theta_i^2 \Vert X_i \Vert^{2}_{4}  = 3 \Vert X_0 \Vert^{2}_{4}  \, .
\]
So, overall,  the following proposition is valid:

\begin{Proposition}
\label{smoothingLemma}  Let $(X_{k})_{k\in {\mathbb{Z}}}$ be a strictly stationary sequence of  martingale differences with finite fourth moment and $G$ a standard normal variable. 
If \eqref{cond-sigma} is satisfied then, there exists a positive constant $c$ such that, for all $T\geq T_{0} \geq 1$ and all $n \geq 1$, 
\begin{equation}
c {\mathbb{E}}\big (\kappa _{\theta }(P_{\langle X,\theta \rangle },P_{G})\big
)\leq \int_{0}^{T_{0}}{\mathbb{E}}\big (|f_{\theta }(t)-\mathrm{e}^{-t^{2}/2}|%
\big )\frac{dt}{t}+\frac{\log (T/T_{0})}{n}+\mathrm{e}^{-T_{0}^{2}/16}+\frac{%
1}{T}\,.  \label{B2}
\end{equation}
\end{Proposition}

If \eqref{condsymmart} is satisfied then \eqref{cond-sigma}
holds and we can apply the smoothing Proposition \ref{smoothingLemma}.
In the rest of the proof we shall take $n \geq 2$ and choose $T=n$ and $T_{0}=4\sqrt{%
\log n}$. 

To derive an upper bound for ${\mathbb{E}}\big (\kappa _{\theta
}(P_{\langle X,\theta \rangle },P_{G})\big )$ of order $1/n$ modulo an
extra-logarithmic term $(\log n)^{2}$, one then needs
to prove that 
\begin{equation}
\int_{0}^{T_{0}}{\mathbb{E}}\big (|f_{\theta }(t)-\mathrm{e}^{-t^{2}/2}|\big
)\frac{dt}{t}\ll \frac{(\log n)^{2 }}{n}\,.  \label{B3mart}
\end{equation}

This will be achieved by a two steps procedure. For any fixed $\theta $,  using Lindeberg's method we
shall replace one by one the variables $\theta_i X_i$ by independent random variables $Y_i(\theta)$ 
taking only two values with some desired characteristics. After that we
shall compare the characteristic function of $\sum_{i=1}^n Y_i(\theta)$ with $\mathrm{e}^{-t^{2}/2}$.

To specify the
distribution of the $Y_i(\theta)$'s, let us mention the following fact.

\begin{Fact}
\label{defYtheta} Let $\sigma^2 >0$ and $\beta_3 \in {\mathbb{R}}$. There
exists a random variable $Y$ taking only 2 values $m$ and $m^{\prime }$
(depending only on $\sigma^2$ and $\beta_3$) and such that ${\mathbb{E}}%
(Y)=0 $, ${\mathbb{E}}(Y^2)=\sigma^2$, ${\mathbb{E}}(Y^3)=\beta_3$ and ${%
\mathbb{E}}(Y^4)= \sigma^4 + \frac{\beta_3^2}{\sigma^2}$.
\end{Fact}

\noindent \textbf{Proof.} According to Lemma 5.1 in \cite{DR08}, we can select $m$ and $m'$ as follows: 
\begin{equation*}
m = \frac{ \beta_3 + \sqrt{\beta_3^2 + 4\sigma^6}}{2 \sigma^2} \, , \,
m^{\prime }=- \frac{\sigma^2}{m} \, , 
\end{equation*}
and consider a r.v. $Y$ with values in $\{m, m'\}$ such that 
\begin{equation*}
{\mathbb{P}} ( Y=m ) =t \text{ and } {\mathbb{P}} ( Y=m^{\prime }) = 1 - t
\, ,
\end{equation*}
where 
\begin{equation*}
t= \frac{ 2 \sigma^6} {4 \sigma^6 + \beta_3 ( \beta_3 + \sqrt{\beta_3^2 + 4
\sigma^6} )} \, .
\end{equation*}

Indeed, in  this case,  by straightforward computations,  ${\mathbb{E}}(Y)=0$, ${\mathbb{E}}(Y^2)=\sigma^2$ and ${\mathbb{E}}%
(Y^3)=\beta_3$. Let us now compute $\Vert Y \Vert_4^4$. Note
that 
\begin{equation*}
m - \sigma^2/m = m+m^{\prime }= \frac{\beta_3}{\sigma^2} \, .
\end{equation*}
Setting $\kappa_3 = \beta_3 / \sigma^2$, we have 
\begin{align*}
\Vert Y \Vert_4^4 & = m^2 ( m^2 t + m^{\prime \, 2} (1-t) ) + m^{\prime \,
2} ( m^{\prime \, 2} - m^2 ) (1-t) \\
& = m^2 \sigma^2 - m^{\prime \, 2} ( m - m^{\prime }) \kappa_3 (1-t) \, .
\end{align*}
But $- m^{\prime \, 2} (1-t) = m^2 t - \sigma^2$ and $m-m^{\prime }= \sqrt{
\kappa_3^2 + 4 \sigma^2}$. Then, simple computations lead to 
\begin{equation*}
\Vert Y \Vert_4^4 = \sigma^4 + \sigma^2 \kappa_3^2 = \sigma^4 + \frac{%
\beta_3^2}{\sigma^2} \, ,
\end{equation*}
which ends the proof of the lemma. \qed

\smallskip

Let $(Y_{i}(\theta ))_{i\geq 1}$ be a sequence of random variables that are
independent for any fixed $\theta $, independent of $(X_{i})_{i\geq 1}$, and
such that, for each $i\geq 1$, the conditional law of $Y_{i}(\theta )$ given 
$\theta $ takes 2 values and is such that ${\mathbb{E}}_{|\theta }(Y_{i}(\theta ))=0$, $%
{\mathbb{E}}_{|\theta }(Y_{i}^{2}(\theta ))=\theta _{i}^{2}{\mathbb{E}}%
(X_{0}^{2})$, ${\mathbb{E}}_{|\theta }(Y_{i}^{3}(\theta ))=\beta
_{i,3}(\theta )$ and ${\mathbb{E}}_{|\theta }(Y_{i}^{4}(\theta ))=\beta
_{i,4}(\theta )$ where 
\begin{equation}
\beta _{k,3}(\theta )=\theta _{k}^{3}{\mathbb{E}}(X_{k}^{3})+3\sum_{\ell
=1}^{k-1}\theta _{\ell }\theta _{k}^{2}{\mathbb{E}}(X_{\ell
}X_{k}^{2})=:\theta _{k}^{3}{\mathbb{E}}(X_{k}^{3})+3{\tilde{\beta}}%
_{k,3}(\theta )  \label{defbeta3k}
\end{equation}%
and 
\begin{equation}
\beta _{k,4}(\theta )=\theta _{k}^{4}({\mathbb{E}}(X_{k}^{2}))^{2}+\frac{%
\beta _{k,3}^{2}(\theta )}{\theta _{k}^{2}{\mathbb{E}}(X_{k}^{2})} = \theta _{k}^{4}+\frac{%
\beta _{k,3}^{2}(\theta )}{\theta _{k}^{2}}\,.
\label{defbeta4k}
\end{equation}%
Note that this is always possible according to Fact \ref{defYtheta} (if $\theta _{k}=0$ we set $\beta _{k,4}(\theta ) =0$ and take $Y_{k}(\theta )=0$).

\medskip

Setting $f_t(x) = \mathrm{e}^{\mathrm{i} t x }$, $M_{k} (\theta) =
\sum_{i=1}^k \theta_i X_i$ and $T_{k,n} (\theta) = \sum_{i=k+1}^n Y_i(
\theta)$, we have 
\begin{multline*}
{\mathbb{E}}_{| \theta} \big ( f_t \big ( \sum_{i=1}^n \theta_i X_i \big ) %
\big ) - {\mathbb{E}}_{| \theta} \big ( f_t \big ( \sum_{i=1}^n Y_i (
\theta) \big ) \big ) \\
= \sum_{k=1}^n \Big \{ {\mathbb{E}}_{| \theta} \big ( f_t \big ( M_{k-1} (
\theta) + \theta_k X_k + T_{k,n} (\theta) \big ) \big ) - {\mathbb{E}}_{|
\theta} \big ( f_t \big ( M_{k-1} ( \theta) + Y_k (\theta) + T_{k,n}
(\theta) \big ) \big ) \Big \} \, .
\end{multline*}
Let 
\begin{equation*}
f_{t,k,n} (x) = {\mathbb{E}}_{| \theta} \big ( f_t \big (x+ T_{k,n} (\theta) %
\big ) \, .
\end{equation*}
This function is in ${\mathcal{C}}^{\infty}$ and all its successive derivatives
are bounded and satisfy: for any $i \geq 0$, $\Vert f^{(i)}_{t,k,n} (x)
\Vert_{\infty} \leq t^i $. By Taylor's expansion and independence between
sequences, it follows that 
\begin{equation}  \label{firstdec}
{\mathbb{E}}_{| \theta} \big ( f_t \big ( \sum_{i=1}^n \theta_i X_i \big ) %
\big ) - {\mathbb{E}}_{| \theta} \big ( f_t \big ( \sum_{i=1}^n Y_i (
\theta) \big ) \big ) = \sum_{i=1}^3 \sum_{k=1}^n I_{i,k} + \sum_{k=1}^n (
R_{1,k}(f_t) - R_{2,k} (f_t)) \, ,
\end{equation}
where the following notations have been used: for any integer $i \geq 1$, 
\begin{equation*}
I_{i,k}= \frac{1}{ i ! } {\mathbb{E}}_{| \theta} \Big ( f^{(i)}_{t,k,n} \big
( M_{k-1} ( \theta) \big ) \big ( \theta^i_k {\mathbb{E}}_{| \theta, k-1} (
X^i_k) - {\mathbb{E}}_{| \theta} ( Y^i_k ( \theta)) \big ) \Big ) \, ,
\end{equation*}
\begin{equation*}
R_{1,k}(f_t) = \frac{1}{6} \int_0^1 (1-s)^3 {\mathbb{E}}_{| \theta} \big \{ %
\theta^4_k X^4_k f^{(4)}_{t,k,n} \big ( M_{k-1} ( \theta) + s \theta_k X_k %
\big ) \big \} ds
\end{equation*}
and 
\begin{equation*}
R_{2,k} (f_t) = \frac{1}{6} \int_0^1 (1-s)^3 {\mathbb{E}}_{| \theta} \big \{ %
Y^4_k (\theta) f^{(4)}_{t,k,n} \big ( M_{k-1} ( \theta) + s Y_k( \theta) %
\big ) \big \} ds \, .
\end{equation*}
Above the notation ${\mathbb{E}}_{| \theta, k-1} (
X^i_k)$ means ${\mathbb{E}} (
X^i_k | \sigma(\theta) \vee {\mathcal F}_{ k-1} )$ (note that by independence between $\theta$ and $X$, ${\mathbb{E}}_{| \theta, k-1} (
X^i_k) = {\mathbb{E}}_{ k-1} (
X^i_k) $). 
Clearly $\sum_{k=1}^n I_{1,k} =0$ (by the martingale property and the fact
that ${\mathbb{E}}_{| \theta}( Y_i( \theta) ) =0$), and 
\begin{equation}  \label{resteTay}
 \sum_{k=1}^n  ( | R_{1,k}(f_t) | + |R_{2,k} (f_t)  |  ) \leq \frac{t^4 }{%
24} \sum_{k=1}^n \big \{ \beta_{k,4} ( \theta) + \theta_k^4 {\mathbb{E}} (
X_k^4) \big \} \, .
\end{equation}
On another hand, since ${\mathbb{E}}_{| \theta} ( Y^2_k( \theta) )
=\theta_k^2 {\mathbb{E}} ( X_0^2)$, 
\begin{equation*}
{\mathbb{E}}_{| \theta} \Big ( f^{(2)}_{t,k,n} (0 ) \big ( \theta^2_k {%
\mathbb{E}}_{| \theta, k-1} ( X^2_k) - {\mathbb{E}}_{| \theta} ( Y^2_k (
\theta)) \big ) \Big ) = \theta_k^2 f^{(2)}_{t,k,n} (0 ) {\mathbb{E}} \big ( 
{\mathbb{E}}_{ k-1} ( X^2_k) - {\mathbb{E}} ( X^2_k) \big ) =0 \, .
\end{equation*}
Hence, since $M_{0} ( \theta) =0$, 
\begin{align}  \label{decI2k}
\sum_{k=1}^n I_{2,k} & = \frac{1}{2} \sum_{k=1}^n \theta^2_k
\sum_{\ell=1}^{k-1} {\mathbb{E}}_{| \theta} \Big \{ \Big ( f^{(2)}_{t,k,n} %
\big ( M_{\ell} ( \theta) \big ) - f^{(2)}_{t,k,n} \big ( M_{\ell-1} (
\theta) \big ) \Big ) \Big ( {\mathbb{E}}_{ \ell} ( X^2_k) - {\mathbb{E}} (
X^2_k ) \Big ) \Big \}  \nonumber \\
& = \frac{1}{2} \sum_{k=1}^n \theta^2_k \sum_{\ell=1}^{k-1} \theta_{\ell} {%
\mathbb{E}}_{| \theta} \Big \{ f^{(3)}_{t,k,n} \big ( M_{\ell-1} ( \theta) %
\big ) X_{\ell} \Big ( {\mathbb{E}}_{ \ell} ( X^2_k) - {\mathbb{E}} ( X^2_k
) \Big )  \Big \}  \nonumber \\
& \quad \quad + \frac{1}{2} \sum_{k=1}^n \theta^2_k \sum_{\ell=1}^{k-1} \theta^2_{\ell}
\int_0^1 (1-s) {\mathbb{E}}_{| \theta} \Big \{ f^{(4)}_{t,k,n} \big ( %
M_{\ell-1} ( \theta) + s \theta_{\ell} X_{\ell}\big ) X^2_{\ell} \Big ( {%
\mathbb{E}}_{ \ell} ( X^2_k) - {\mathbb{E}} ( X^2_k ) \Big )  \Big \} \nonumber
\\
& := \sum_{k=1}^n I_{2,3,k} + \sum_{k=1}^n I_{2,4,k} \, .
\end{align}
We have 
\begin{equation}  \label{borneI24}
\Big | \sum_{k=1}^n I_{2,4,k} \Big | \leq \frac{t^4}{4} \sum_{k=1}^n
\theta^2_k \sum_{\ell=1}^{k-1} \theta^2_{\ell} \gamma_{2,2} (k-\ell) \leq 
\frac{t^4}{4} \sum_{k=1}^n \theta^4_k \sum_{v=1}^{n} \gamma_{2,2} (v) \, .
\end{equation}
Next we deal with the quantity $\sum_{k=1}^n \big ( I_{3,k} + I_{2,3,k} \big )$. Since ${\mathbb{E}}_{| \theta} ( Y^3_i( \theta) ) =\beta_{i,3} (
\theta)$, by the definition of ${{\beta}}%
_{k,3}(\theta ) $ and ${\tilde{\beta}}%
_{k,3}(\theta ) $ in \eqref{defbeta3k}, we have 
\begin{multline*}
\sum_{k=1}^n \big ( I_{3,k} + I_{2,3,k} \big )=  \frac{1}{ 6 } \sum_{k=1}^n   \theta^3_k   {\mathbb{E}}_{| \theta} \Big ( f^{(3)}_{t,k,n} \big
( M_{k-1} ( \theta) \big ) \big ( {\mathbb{E}}_{ k-1} (
X^3_k) - {\mathbb{E}}( X^3_k) \big ) \Big ) \\
+ \frac{1}{2} \sum_{k=1}^n \theta^2_k \sum_{\ell=1}^{k-1} \theta_{\ell} {%
\mathbb{E}}_{| \theta} \Big \{ f^{(3)}_{t,k,n} \big ( M_{\ell-1} ( \theta) %
\big ) X_{\ell} \Big ( {\mathbb{E}}_{ \ell} ( X^2_k) - {\mathbb{E}} ( X^2_k
) \Big )  \Big \}  \\  -  \frac{1}{2} \sum_{k=1}^n {\tilde \beta}_{k,3} ( \theta) {\mathbb{E}}_{|
\theta} \Big ( f^{(3)}_{t,k,n} \big ( M_{k-1} ( \theta) \big ) \Big ) \, .
\end{multline*}
Next define $Z_{\ell, k} = X_{\ell} \big ( {\mathbb{E}}_{ \ell} ( X^2_k) - {\mathbb{E}} (
X^2_k ) \big ) $ and 
\begin{align*}
J_n& := \frac{1}{2} \sum_{k=1}^n \Big \{ {\tilde \beta}_{k,3} ( \theta) {\mathbb{E}}_{|
\theta} \Big ( f^{(3)}_{t,k,n} \big ( M_{k-1} ( \theta) \big ) \Big ) -
\sum_{\ell=1}^{k-1} \theta_k^2 \theta_{\ell} {\mathbb{E}}_{| \theta} \Big ( %
f^{(3)}_{t,k,n} \big ( M_{\ell-1} ( \theta) \big ) \Big ) {\mathbb{E}}  ( Z_{\ell, k} ) \Big \} \\
& =  \frac{1}{2} \sum_{k=1}^n \sum_{\ell=1}^{k-1} \theta_k^2 \theta_{\ell} \Big \{ {\mathbb{%
E}}_{| \theta} \Big ( f^{(3)}_{t,k,n} \big ( M_{k-1} ( \theta) \big ) \Big ) %
- {\mathbb{E}}_{| \theta} \Big ( f^{(3)}_{t,k,n} \big ( M_{\ell-1} ( \theta) %
\big ) \Big ) \Big \} {\mathbb{E}} ( X_{\ell} X^2_k) \, .
\end{align*}
For $\ell \leq k$, 
\begin{align*}
& {\mathbb{E}}_{| \theta} \Big ( f^{(3)}_{t,k,n} \big ( M_{k-1} ( \theta) \big
) \Big ) - {\mathbb{E}}_{| \theta} \Big ( f^{(3)}_{t,k,n} \big ( M_{\ell-1}
( \theta) \big ) \Big )  \\ & = \sum_{u=\ell}^{ k-1} \Big \{ {\mathbb{E}}_{|
\theta} \Big ( f^{(3)}_{t,k,n} \big ( M_{u} ( \theta) \big ) \Big ) - {%
\mathbb{E}}_{| \theta} \Big ( f^{(3)}_{t,k,n} \big ( M_{u-1} ( \theta) \big
) \Big ) \Big \} \\
& = \sum_{u=\ell}^{ k-1} \theta_u {\mathbb{E}}_{| \theta} \Big ( %
f^{(4)}_{t,k,n} \big ( M_{u-1} ( \theta) \big ) \Big ) {\mathbb{E%
}}_{u-1} (X_u)  + \int_0^1 (1-s) \sum_{u=\ell}^{ k-1} \theta^2_u {%
\mathbb{E}}_{| \theta} \Big ( f^{(5)}_{t,k,n} \big ( M_{u-1} ( \theta)+s X_u %
\big ) X_u^2 \Big ) ds \, .
\end{align*}
By the Martingale property, ${\mathbb{E%
}}_{u-1} (X_u) =0$. Therefore
\begin{equation}  \label{BorneJn}
|J_n| \leq \frac{ t^5}{4} \sum_{k=1}^n \sum_{\ell=1}^{k-1} \sum_{u=\ell}^{
k-1} \theta_k^2 |\theta_{\ell} | \theta^2_u \gamma_{1,2} ( k - \ell ) \, .
\end{equation}
Next
\begin{multline*}
\sum_{k=1}^n \big ( I_{3,k} + I_{2,3,k} \big ) +J_n =  \frac{1}{ 6 } \sum_{k=1}^n   \theta^3_k   {\mathbb{E}}_{| \theta} \Big ( f^{(3)}_{t,k,n} \big
( M_{k-1} ( \theta) \big ) \big ( {\mathbb{E}}_{ k-1} (
X^3_k) - {\mathbb{E}}( X^3_k) \big ) \Big ) \\
+ \frac{1}{2} \sum_{k=1}^n \theta^2_k \sum_{\ell=1}^{k-1} \theta_{\ell} {%
\mathbb{E}}_{| \theta} \Big \{ f^{(3)}_{t,k,n} \big ( M_{\ell-1} ( \theta) %
\big ) X_{\ell} \Big ( {\mathbb{E}}_{ \ell} ( X^2_k) - {\mathbb{E}} ( X^2_k
) \Big )  \Big \}  \\
  -
\sum_{\ell=1}^{k-1} \theta_k^2 \theta_{\ell} {\mathbb{E}}_{| \theta} \Big ( %
f^{(3)}_{t,k,n} \big ( M_{\ell-1} ( \theta) \big ) \Big ) {\mathbb{E}} \Big
( X_{\ell} \big ( {\mathbb{E}}_{ \ell} ( X^2_k) - {\mathbb{E}} ( X^2_k ) %
\big ) \Big ) \Big \}
\, .
\end{multline*}
Therefore, recalling the notation $Z_{\ell, k} = X_{\ell} \big ( {\mathbb{E}}_{ \ell} ( X^2_k) - {\mathbb{E}} (
X^2_k ) \big ) $, 
and using that, by the martingale property,  $ {\mathbb{E}}_{ k-1}  (Z_{k, k} )=  {\mathbb{E}}_{ k-1} (
X^3_k) $,  we get
\begin{equation*}
\sum_{k=1}^n \big ( I_{3,k} + I_{2,3,k} \big ) + J_n = \frac{1}{6} \sum_{k=1}^n \theta^2_k \sum_{\ell=1}^{k} \theta_{\ell} ( 1+ 2 \mathbf{1}_{k
\neq \ell})  {\mathbb{E}}_{|
\theta} \Big \{ f^{(3)}_{t,k,n} \big ( M_{\ell-1} ( \theta) \big ) \big ( %
Z_{\ell,k} - {\mathbb{E}} ( Z_{\ell,k} ) \big ) \Big \} \, .
\end{equation*}
Hence, since ${\mathbb{E}}_{| \theta} \Big \{    f^{(3)}_{t,k,n} \big ( 0 %
\big )  \big (  Z_{\ell,k} - {\mathbb{E}} ( Z_{\ell,k} ) \big )  \Big \} =0$%
, 
\begin{align*}
& \sum_{k=1}^n \big ( I_{3,k} + I_{2,3,k} \big ) + J_n \\
& = \frac{1}{6} \sum_{k=1}^n \theta^2_k \sum_{\ell=1}^{k} ( 1+ 2 \mathbf{1}_{k
\neq \ell})\sum_{u=1}^{\ell - 1}\theta_{\ell} {\mathbb{E}}_{| \theta} \Big
\{ \Big ( f^{(3)}_{t,k,n} \big ( M_{\ell-u} ( \theta) \big ) -
f^{(3)}_{t,k,n} \big ( M_{\ell-u-1} ( \theta) \big ) \Big ) \big ( %
Z_{\ell,k} - {\mathbb{E}} ( Z_{\ell,k} ) \big ) \Big \} \\
& = \frac{1}{6} \sum_{k=1}^n \theta^2_k \sum_{\ell=1}^{k} ( 1+ 2 \mathbf{1}_{k
\neq \ell}) \sum_{u=1}^{\ell - 1}\theta_{\ell} \theta_{\ell - u} \int_0^1 L_{k,\ell, u} (\theta , s)ds \, .
\end{align*}
where
\[
L_{k,\ell, u} (\theta , s) := {\mathbb{E}}_{| \theta} \Big \{ f^{(4)}_{t,k,n} \big ( M_{\ell-u-1} ( \theta)
+ s \theta_{\ell - u} X_{\ell - u} \big ) X_{\ell - u} \big ( Z_{\ell,k} - {%
\mathbb{E}} ( Z_{\ell,k} ) \big ) \Big \} \, .
\]
Let $u \geq 1$. By the Martingale property,  ${\mathbb{E}}
_{\ell - u } (Z_{\ell, k})  = {\mathbb{E}}
_{\ell - u }  (  X_{\ell} X^2_k) $, which implies that 
\[
\Big | L_{k,\ell, u} (\theta , s)\Big | 
\leq t^4 \Vert X_{\ell - u} \big ( {\mathbb{E}} _{\ell - u } (
X_{\ell}X^2_k) - {\mathbb{E}} ( X_{\ell}X^2_k) \big ) \Vert_1 \leq t^4
\gamma_{1,3} (u) \, .
\]
On another hand 
\begin{multline*}
\Big | L_{k,\ell, u} (\theta , s)\Big | 
\leq \Big | {\mathbb{E}}_{| \theta} \Big \{ f^{(4)}_{t,k,n} \big ( %
M_{\ell-u-1} ( \theta) + s \theta_{\ell - u} X_{\ell - u} \big ) X_{\ell -
u} X_{\ell } \big ( {\mathbb{E}}_{ \ell} ( X^2_k) - {\mathbb{E}} ( X^2_k ) %
\big ) \Big \} \Big | \\
+ \Big | {\mathbb{E}}_{| \theta} \Big \{ f^{(4)}_{t,k,n} \big ( M_{\ell-u-1}
( \theta) + s \theta_{\ell - u} X_{\ell - u} \big ) X_{\ell - u} \Big \} %
\Big | \Big \Vert X_{\ell } \big ( {\mathbb{E}}_{ \ell} ( X^2_k) - {\mathbb{E%
}} ( X^2_k ) \big ) \Big \Vert_1 \\
\leq t^4 \Vert X_{0} X_{u } \big ( {\mathbb{E}}_u ( X^2_{k-\ell +u} - {%
\mathbb{E}} ( X^2_0) \big ) \Vert_1 + t^4 \Vert X_0 \Vert_1 \Vert X_{0} \big
( {\mathbb{E}}_{ 0} ( X^2_{k- \ell}) - {\mathbb{E}} ( X^2_0 ) \big ) \Vert_1 \\
\leq t^4 ( \gamma_{2,2} (k- \ell ) + \gamma_{1,2} (k - \ell ) ) \, .
\end{multline*}
Bearing in mind Definition \eqref{defgamma}, we get 
\begin{equation}  \label{borneI3+I23}
\Big | \sum_{k=1}^n \big ( I_{3,k} + I_{2,3,k} \big ) + J_n \Big | \leq t^4
\sum_{k=1}^n \theta^2_k \sum_{\ell=1}^{k} \sum_{u=1}^{\ell - 1} |
\theta_{\ell} \theta_{\ell - u} | \Big ( \gamma (u) \wedge \gamma (k - \ell
) \Big ) \, .
\end{equation}
Starting from \eqref{firstdec} and taking into account \eqref{resteTay}, %
\eqref{decI2k}, \eqref{borneI24}, \eqref{BorneJn} and \eqref{borneI3+I23},
we derive that 
\begin{multline}  \label{concludec}
\Big | {\mathbb{E}}_{| \theta} \big ( f_t \big ( \sum_{i=1}^n \theta_i X_i %
\big ) \big ) - {\mathbb{E}}_{| \theta} \big ( f_t \big ( \sum_{i=1}^n Y_i (
\theta) \big ) \big ) \Big | \leq \frac{t^4 }{12} \sum_{k=1}^n \big \{ %
\beta_{k,4} ( \theta) + \theta_k^4 {\mathbb{E}} ( X_0^4) \big \}  \\ + \frac{t^4%
}{4} \sum_{k=1}^n \theta^4_k \sum_{v=1}^{n} \gamma (v) 
+ t^4 \sum_{k=1}^n \theta^2_k \sum_{\ell=1}^{k} \sum_{u=1}^{\ell - 1} |
\theta_{\ell} \theta_{\ell - u} | \Big ( \gamma (u) \wedge \gamma (k - \ell
) \Big )  \\ + \frac{ t^5}{4} \sum_{k=1}^n \sum_{\ell=1}^{k-1} \sum_{u=\ell}^{
k-1} \theta_k^2 |\theta_{\ell} | \theta^2_u \gamma( k - \ell ) \, .
\end{multline}
By Young's inequality, 
\begin{equation*}
\theta_k^2| \theta_{\ell} \theta_{\ell - u} | \leq \frac{1}{2 \sqrt{2}} \big
( \theta_k^4 + \theta^4_{\ell} + \theta^4_{\ell - u} \big ) \ \text{ and } \
\theta_k^2 |\theta_{\ell} | \theta^2_u \leq \frac{2}{5} \big ( | \theta_k|^5
+ | \theta_u|^5 + | \theta_\ell |^5 \big ) \, .
\end{equation*}
Now, for any $m \geq 1$, ${\mathbb{E}} ( |\theta_v|^m) \leq c_m n^{-m/2}$.
Hence 
\begin{equation*}
{\mathbb{E}} ( \theta_k^2| \theta_{\ell} \theta_{\ell - u} | ) \ll n^{-2} 
\text{ and } {\mathbb{E}} ( \theta_k^2 |\theta_{\ell} | \theta^2_u ) \ll
n^{-5/2} \, .
\end{equation*}
In addition, 
\begin{equation}  \label{Bornebeta4}
{\mathbb{E}} ( \beta_{k,4} ( \theta) ) \ll n^{-2} \Big ( 1 + \sum_{\ell
=1}^{k-1} \gamma ( \ell ) \Big )^2 \, .
\end{equation}
Indeed, 
\begin{equation*}
\beta_{k,4} ( \theta) \leq \theta_k^4 ( 1 + 2 \Vert X_0 \Vert_3^2) + 18
\theta^2_{k} \Big ( \sum_{\ell =1}^{k-1} \theta_{\ell} {\mathbb{E}} (
X_{\ell} X_k^2) \Big )^2 \, .
\end{equation*}
Then we use the fact that $\sup_{1 \leq v \leq n} {\mathbb{E}} ( \theta_v^4)
\leq C n^{-2}$ and that 
\begin{equation*}
{\mathbb{E}} \Big ( \theta^2_{k} \Big ( \sum_{\ell =1}^{k-1} \theta_{\ell} {%
\mathbb{E}} ( X_{\ell} X_k^2) \Big )^2 \Big ) \leq \sum_{\ell , \ell^{\prime
}=1}^{k-1} \gamma ( k - \ell ) \gamma ( k - \ell^{\prime }) {\mathbb{E}}
(|\theta^2_{k} \theta_{\ell} \theta_{\ell^{\prime }} |) \leq C n^{-2} \Big ( %
\sum_{\ell =1}^{k-1} \gamma ( \ell ) \Big )^2 \, .
\end{equation*}

So, overall, starting from \eqref{concludec}, we derive that 
\begin{multline*}
{\mathbb{E}}\Big |{\mathbb{E}}_{|\theta }\big (f_{t}\big (%
\sum_{i=1}^{n}\theta _{i}X_{i}\big )\big )-{\mathbb{E}}_{|\theta }\big (f_{t}%
\big (\sum_{i=1}^{n}Y_{i}(\theta )\big )\big )\Big |\ll t^4 n^{-1}\Big (%
1+\sum_{\ell =1}^{n}\gamma (\ell )\Big )^{2} \\
+t^{4}n^{-1}\sum_{v=1}^{n}\sum_{u=1}^{n}\gamma (u)\wedge \gamma
(v)+t^{5}n^{-3/2}\sum_{v=1}^{n}v\gamma (v)\,.
\end{multline*}%
Hence, the following bound is valid: 
\begin{multline*}
{\mathbb{E}}\Big |{\mathbb{E}}_{|\theta }\big (f_{t}\big (%
\sum_{i=1}^{n}\theta _{i}X_{i}\big )\big )-{\mathbb{E}}_{|\theta }\big (f_{t}%
\big (\sum_{i=1}^{n}Y_{i}(\theta )\big )\big )\Big | \\ \ll \frac{t^{4}}{n}\Big
\{1+\Big (\sum_{\ell =1}^{n}\gamma (\ell )\Big )^{2}+\sum_{v=1}^{n}v\gamma
(v)+\frac{t}{\sqrt{n}}\sum_{v=1}^{n}v\gamma (v)\Big \}\,,
\end{multline*}%
implying that 
\begin{equation*}
\int_{0}^{T_{0}}t^{-1}{\mathbb{E}}\Big |{\mathbb{E}}_{|\theta }\big (f_{t}%
\big (\sum_{i=1}^{n}\theta _{i}X_{i}\big )\big )-{\mathbb{E}}_{|\theta }\big
(f_{t}\big (\sum_{i=1}^{n}Y_{i}(\theta )\big )\big )\Big |dt\ll \frac{%
T_{0}^{4}}{n}\Big \{1+\Big (\sum_{\ell =1}^{n}\gamma (\ell )\Big )%
^{2}+\sum_{v=1}^{n}v\gamma (v)\Big \}\,.
\end{equation*}%
Since $\sum_{v\geq 1}v\gamma (v)<\infty $, it follows that 
\begin{equation}
\int_{0}^{T_{0}}{\mathbb{E}}\Big ( \Big |{\mathbb{E}}_{|\theta }\big (\mathrm{e}^{%
\mathrm{i}t\sum_{k=1}^{n}X_{k}(\theta )}\big )-{\mathbb{E}}_{|\theta }\big (%
\mathrm{e}^{\mathrm{i}t\sum_{k=1}^{n}Y_{k}(\theta )}\big ) \Big | \Big )\frac{dt}{t}%
\ll \frac{(\log n)^{2}}{n}\,.  \label{conclusion1formart}
\end{equation}%
Therefore the upper bound \eqref{B3mart} will follow from %
\eqref{conclusion1formart} provided one can prove that 
\begin{equation}
\int_{0}^{T_{0}}{\mathbb{E}}\Big (\Big |{\mathbb{E}}_{|\theta }\big (\mathrm{%
e}^{\mathrm{i}t\sum_{k=1}^{n}Y_{k}(\theta )}\big )-\mathrm{e}^{-t^{2}/2}\Big
|\Big )\frac{dt}{t}\ll \frac{(\log n)^{2}}{n}\,.  \label{B3indepart}
\end{equation}%

With this aim, we shall adapt the proof of \cite[Lemma 2.1]{KS}. Their result cannot be applied directly since $Y_k(\theta)$ is not of the form $\theta_k \eta_k$ where $(\eta_k)_{k \in \mathbb Z}$ is a sequence of independent r.v.'s independent of $\theta$. Let 
\begin{equation*}
\Gamma _{n}(\theta )=\Big \{\max_{1\leq k\leq n}|\beta _{k,3}(\theta
)|T_{0}\leq 1\Big \}\cap \Big \{T_{0}^{3}\Big |\sum_{k=1}^{n}\beta
_{k,3}(\theta )\Big |\leq 1\Big \}\cap \Big \{T_{0}^{4}\sum_{k=1}^{n}\beta
_{k,4}(\theta )\leq 1\Big \} \, .
\end{equation*}%
Since, when $\theta $ is fixed but in ${S^{n-1}}$, $(Y_{k}(\theta ))_{1\leq k\leq n}$ are
independent random variables that are centered, in ${\mathbb{L}}^{4}$ and
such that $\sum_{k=1}^{n}{\mathbb{E}}(Y_{k}^{2}(\theta ))=1$, we infer that,
by standard arguments (see  the proof of \cite[Lemma 2.1]{KS}),
the following estimate holds : for any positive $t$ such that $t\leq T_{0}$, 
\begin{equation}
\Big |{\mathbb{E}}_{|\theta }\Big (\mathrm{e}^{\mathrm{i}%
t\sum_{k=1}^{n}Y_{k}(\theta )}\Big )-\mathrm{e}^{-t^{2}/2}\Big |\mathbf{1}%
_{\Gamma _{n}(\theta )}\ll \mathrm{e}^{-t^{2}/2}\Big (t^{3}\Big |%
\sum_{k=1}^{n}\beta _{k,3}(\theta )\Big |+t^{4}\sum_{k=1}^{n}\beta
_{k,4}(\theta )\Big )\,.  \label{B3indepartP1}
\end{equation}%
Let $a(u)={\mathbb{E}}(X_{0}X_{u}^{2})$. Note that 
\begin{equation*}
{\mathbb{E}}\Big |\sum_{k=1}^{n}\beta _{k,3}(\theta )\Big |\leq \Vert
X_{0}\Vert _{3}^{3}{\mathbb{E}}\Big |\sum_{k=1}^{n}\theta _{k}^{3}\Big |+ 3 {%
\mathbb{E}}\Big |\sum_{k=1}^{n}\theta _{k}^{2}\sum_{\ell =1}^{k-1}a(k-\ell
)\theta _{\ell }\Big |\,.
\end{equation*}%
We shall use the fact that $(\theta _{1},\ldots ,\theta _{n})=^{\mathcal{D}%
}(\xi _{1},\ldots ,\xi _{n})/\Vert \xi \Vert_e $ where $(\xi _{i})_{i\geq 1}$
is a sequence of iid ${\mathcal{N}}(0,1)$-distributed r.v.'s and $\Vert \xi \Vert_e
^{2}=\sum_{i=1}^{n}\xi _{i}^{2}$. Note first that, since there exists $K>0$ such that for any $n >6$,   ${\mathbb{E}}\Big(%
n^{3}/\Vert \xi \Vert_e ^{6}\Big )\leq K$ (by the properties of the $\chi^2$-distribution),  
\begin{equation*}
{\mathbb{E}}\Big |\sum_{k=1}^{n}\theta _{k}^{3}\Big |={\mathbb{E}}\Big |%
\frac{1}{\Vert \xi \Vert_e ^{3}}\sum_{k=1}^{n}\xi _{k}^{3}\Big |\leq   {\mathbf 1}_{n\leq 6}  + {\mathbf 1}_{n>6} {\Big ({%
\mathbb{E}}\Big(\frac{n^{3}}{\Vert \xi \Vert_e ^{6}}\Big )\Big )^{1/2}}\Big (%
\frac{1}{n^{3}}{\mathbb{E}}\Big (\sum_{k=1}^{n}\xi _{k}^{3}\Big )^{2}\Big )%
^{1/2}\ll n^{-1}\,.
\end{equation*}%
Next 
\begin{multline*}
{\mathbb{E}}\Big |\sum_{k=1}^{n}\theta _{k}^{2}\sum_{\ell =1}^{k-1}a(k-\ell
)\theta _{\ell }\Big |={\mathbb{E}}\Big |\frac{1}{\Vert \xi \Vert_e ^{3}}%
\sum_{k=1}^{n}\xi _{k}^{2}\sum_{\ell =1}^{k-1}a(k-\ell )\xi _{\ell }\Big | \\
\ll  {\mathbf 1}_{n\leq 6}  \sum_{\ell=1}^5 |a(\ell) |  + {\mathbf 1}_{n>6}{\Big ({\mathbb{E}}\Big(\frac{n^{3}}{\Vert \xi \Vert_e ^{6}}\Big )\Big )%
^{1/2}}\Big (\frac{1}{n^{3}}{\mathbb{E}}\Big (\sum_{k=1}^{n}\xi
_{k}^{2}\sum_{\ell =1}^{k-1}a(k-\ell )\xi _{\ell }\Big )^{2}\Big )^{1/2}\,.
\end{multline*}%
Now, by independence, 
\begin{multline*}
{\mathbb{E}}\Big (\sum_{k=1}^{n}\xi _{k}^{2}\sum_{\ell =1}^{k-1}a(k-\ell
)\xi _{\ell }\Big )^{2}=\sum_{k=1}^{n}{\mathbb{E}}(\xi _{k}^{4})\sum_{\ell
=1}^{k-1}a^{2}(k-\ell ){\mathbb{E}}(\xi _{\ell }^{2}) \\
+2\sum_{k=1}^{n-1}\sum_{k^{\prime }=k+1}^{n}{\mathbb{E}}(\xi _{k^{\prime
}}^{2})\sum_{\ell =1}^{k-1}\sum_{\ell ^{\prime }=1}^{k^{\prime }-1}a(k-\ell
)a(k^{\prime }-\ell ^{\prime }){\mathbb{E}}(\xi _{\ell }\xi _{k}^{2}\xi
_{\ell ^{\prime }})\,.
\end{multline*}%
By independence again, 
\begin{equation*}
\sum_{k=1}^{n-1}\sum_{k^{\prime }=k+1}^{n}\sum_{\ell =1}^{k-1}\sum_{\ell
^{\prime }=1}^{k^{\prime }-1}a(k-\ell )a(k^{\prime }-\ell ^{\prime }){%
\mathbb{E}}(\xi _{\ell }\xi _{k}^{2}\xi _{\ell ^{\prime
}})=\sum_{k=1}^{n-1}\sum_{k^{\prime }=k+1}^{n}\sum_{\ell =1}^{k-1}a(k-\ell
)a(k^{\prime }-\ell )\,.
\end{equation*}%
Hence 
\begin{equation*}
{\mathbb{E}}\Big (\sum_{k=1}^{n}\xi _{k}^{2}\sum_{\ell =1}^{k-1}a(k-\ell
)\xi _{\ell }\Big )^{2}\ll n\Big (\sum_{\ell =1}^{n}\gamma (\ell )\Big )%
^{2}\,.
\end{equation*}%
So overall, it follows that 
\begin{equation}
{\mathbb{E}}\Big |\sum_{k=1}^{n}\beta _{k,3}(\theta )\Big |\ll n^{-1}\Big (%
1+\sum_{\ell =1}^{n}\gamma (\ell )\Big )\,.  \label{boundsumbetak3}
\end{equation}%
Starting from \eqref{B3indepartP1} and taking into account %
\eqref{boundsumbetak3} and \eqref{Bornebeta4}, we derive that 
\begin{equation}
\int_{0}^{T_{0}}{\mathbb{E}}\Big (\Big |{\mathbb{E}}_{|\theta }\Big (\mathrm{%
e}^{\mathrm{i}t\sum_{k=1}^{n}Y_{k}(\theta )}\Big )-\mathrm{e}^{-t^{2}/2}\Big
|\mathbf{1}_{\Gamma _{n}(\theta )}\Big )\frac{dt}{t}\ll n^{-1}\Big (%
1+\sum_{\ell =1}^{n}\gamma (\ell )\Big )^{2}\,.  \label{B3indepartP2}
\end{equation}%
Next, note that 
%
\begin{multline}
{\mathbb{P}}\big (\Gamma _{n}^{c}(\theta )\big )\leq \sum_{k=1}^{n}{\mathbb{P%
}}\big (T_{0}|\beta _{k,3}(\theta )|>1\big )+{\mathbb{P}}\Big (T_{0}^{3}\Big
|\sum_{k=1}^{n}\beta _{k,3}(\theta )\Big |>1\Big ) \\+{\mathbb{P}}\Big (%
T_{0}^{4}\Big |\sum_{k=1}^{n}\beta _{k,4}(\theta )\Big |>1\Big )\,.
\label{bornePgammanc}
\end{multline}%
We first deal with the first term in the right-hand side of \eqref{bornePgammanc}. By Markov's
inequality, 
\begin{equation*}
\sum_{k=1}^{n}{\mathbb{P}}\big (T_{0}|\beta _{k,3}(\theta )|>1\big )\ll
T_{0}^{2}\sum_{k=1}^{n}\Big ({\mathbb{E}}(\theta _{k}^{6})+{\mathbb{E}}\Big
\{\theta _{k}^{4}\Big (\sum_{\ell =1}^{k-1}\theta _{\ell }a(k-\ell )\Big )%
^{2}\Big \}\Big )\,.
\end{equation*}%
Now $ \sup_{1 \leq k \leq n}{\mathbb{E}}(\theta _{k}^{6})\ll n^{-3}$, and 
\begin{equation*}
{\mathbb{E}}\Big \{\theta _{k}^{4}\Big (\sum_{\ell =1}^{k-1}\theta _{\ell
}a(k-\ell )\Big )^{2}\Big \}\leq \sum_{\ell ,\ell ^{\prime
}=1}^{k-1}|a(k-\ell )||a(k-\ell ^{\prime })|{\mathbb{E}}(\theta
_{k}^{4}|\theta _{\ell }\theta _{\ell ^{\prime }}|)\ll n^{-3}\Big (%
\sum_{\ell =1}^{k-1}\gamma (\ell )\Big )^{2}\,.
\end{equation*}%
Hence 
\begin{equation}
\sum_{k=1}^{n}{\mathbb{P}}\big (T_{0}|\beta _{k,3}(\theta )|>1\big )\ll
T_{0}^{2}n^{-2}\Big (1+\sum_{\ell =1}^{n}\gamma (\ell )\Big )^{2}\,.
\label{bornegammanc1}
\end{equation}%
Starting from \eqref{bornePgammanc}, applying Markov's inequality and taking
into account the upper bounds \eqref{bornegammanc1}, \eqref{boundsumbetak3}
and \eqref{Bornebeta4}, we derive that 
\begin{equation}
{\mathbb{P}}\big (\Gamma _{n}^{c}(\theta )\big )\ll T_{0}^{4}n^{-1}\Big (%
1+\sum_{\ell =1}^{n}\gamma (\ell )\Big )^{2}\,.  \label{bornePgammanc2}
\end{equation}%
On another hand, note that 
\begin{multline*}
\Big |{\mathbb{E}}_{|\theta }\Big (\mathrm{e}^{\mathrm{i}%
t\sum_{k=1}^{n}Y_{k}(\theta )}\Big )-\mathrm{e}^{-t^{2}/2}\Big | \\ =\Big |%
\prod_{k=1}^{n}{\mathbb{E}}_{|\theta }\Big (\mathrm{e}^{\mathrm{i}%
tY_{k}(\theta )}\Big )-\prod_{k=1}^{n}\mathrm{e}^{-t^{2}\theta _{k}^{2}/2}%
\Big | 
\leq \sum_{k=1}^{n}\Big |{\mathbb{E}}_{|\theta }\Big (\mathrm{e}^{\mathrm{i}%
tY_{k}(\theta )}\Big )-\mathrm{e}^{-t^{2}\theta _{k}^{2}/2}\Big |\,.
\end{multline*}
But $ \big |  \mathrm{e}^{-t^{2}\theta _{k}^{2}/2}  -1 + t^{2}\theta _{k}^{2}/2  \big |  \leq  t^{4}\theta _{k}^{4}/8 $ and 
\[
\Big | {\mathbb{E}}_{|\theta }\Big (\mathrm{e}^{\mathrm{i}%
tY_{k}(\theta )}\Big ) -1 -it   {\mathbb{E}}_{|\theta } (Y_{k}(\theta ) )  + \frac{t^2}{2}  {\mathbb{E}}_{|\theta } (Y^2_{k}(\theta ) )  \Big |  \leq  \frac{|t|^3}{6}  {\mathbb{E}}_{|\theta } ( | Y_{k}(\theta )  |^3)  \, .
\]
Therefore, since $ {\mathbb{E}}_{|\theta } (Y_{k}(\theta ) ) =0$ and $ {\mathbb{E}}_{|\theta } (Y^2_{k}(\theta ) ) = \theta_k^2$, 
\[
\Big |{\mathbb{E}}_{|\theta }\Big (\mathrm{e}^{\mathrm{i}%
t\sum_{k=1}^{n}Y_{k}(\theta )}\Big )-\mathrm{e}^{-t^{2}/2}\Big | \leq \frac{%
|t|^{3}}{6}\sum_{k=1}^{n}{\mathbb{E}}_{|\theta }(|Y_{k}^{3}(\theta )|)+\frac{%
t^{4}}{8}\sum_{k=1}^{n}\theta _{k}^{4}\,.
\]
Now, by H\"older's inequality,  ${\mathbb{E}}_{|\theta }(|Y_{k}^{3}(\theta )|)\leq  \Big (  {\mathbb{E}}_{|\theta }(|Y_{k}^{4}(\theta )|)  \Big )^{3/4}= \beta
_{k,4}^{3/4}(\theta )$. On another hand, by the definition of $ \beta
_{k,4} (\theta )$, using the fact that  $|\theta _{k}|\leq 1$, we have   $\theta
_{k}^{4}\leq  |\theta
_{k}^{3}|  \leq \beta _{k,4}^{3/4}(\theta )$. Therefore 
\begin{equation}
\big |{\mathbb{E}}_{|\theta }\big (\mathrm{e}^{\mathrm{i}%
t\sum_{k=1}^{n}Y_{k}(\theta )}\big )-\mathrm{e}^{-t^{2}/2}\big |\leq \frac{%
|t|^{3}}{6}(1+|t|)\sum_{k=1}^{n}\beta _{k,4}^{3/4}(\theta )\,.
\label{bornePgammanc2bis}
\end{equation}%
Taking into account \eqref{bornePgammanc2}, \eqref{bornePgammanc2bis} and
the fact that $\sum_{\ell \geq 1}\gamma (\ell )<\infty $, we infer that for
any $r\geq 1$, 
\begin{equation}
\int_{0}^{T_{0}}{\mathbb{E}}\Big (\Big |{\mathbb{E}}_{|\theta }\Big (\mathrm{%
e}^{\mathrm{i}t\sum_{k=1}^{n}Y_{k}(\theta )}\Big )-\mathrm{e}^{-t^{2}/2}\Big
|\mathbf{1}_{\Gamma _{n}^{c}(\theta )}\Big )\frac{dt}{t}
\label{interbornePgammanc2} \\
\ll T_{0}^{4}\sum_{k=1}^{n}\Big ({\mathbb{E}}\big (\beta
_{k,4}^{3r/4}(\theta )\big )\Big )^{\frac{1}{r}}\Big (T_{0}^{4}n^{-1}\Big )%
^{\frac{r-1}{r}}\,.
\end{equation}%
Proceeding as to prove \eqref{Bornebeta4}, we infer that if $\sum_{\ell
\geq 1}\gamma (\ell )<\infty $, for any $m\geq 1$,  there exists a finite constant $C_m$ such that for any $n \geq 1$, 
\begin{equation}
\max_{1\leq k\leq n}{\mathbb{E}}(\beta _{k,4}^{m}(\theta ))\leq C_m n^{-2m}\,.
\label{ibornebeta4powerm}
\end{equation}%
So, for any $r\geq 4/3$, 
\begin{equation*}
\int_{0}^{T_{0}}{\mathbb{E}}\Big (\Big |{\mathbb{E}}_{|\theta }\Big (\mathrm{%
e}^{\mathrm{i}t\sum_{k=1}^{n}Y_{k}(\theta )}\Big )-\mathrm{e}^{-t^{2}/2}\Big
|\mathbf{1}_{\Gamma _{n}^{c}(\theta )}\Big )\frac{dt}{t}\ll T_{0}^{8-4/r}%
\frac{n^{1/r}}{n^{3/2}}\,.
\end{equation*}%
Taking $r>2$ in the inequality above, we derive that 
\begin{equation*}
\int_{0}^{T_{0}}{\mathbb{E}}\Big (\Big |{\mathbb{E}}_{|\theta }\Big (\mathrm{%
e}^{\mathrm{i}t\sum_{k=1}^{n}Y_{k}(\theta )}\Big )-\mathrm{e}^{-t^{2}/2}\Big
|\mathbf{1}_{\Gamma _{n}^{c}(\theta )}\Big )\frac{dt}{t}\ll \frac{1}{n}\,,
\end{equation*}%
which combined with \eqref{B3indepartP2} implies \eqref{B3indepart} in case $%
\sum_{\ell \geq 1}\gamma (\ell )<\infty $. This ends the proof of the upper bound \eqref{concluBEmartthm}. \qed

\subsection{Proof of the upper bound \eqref{concluBEmartthm2}} 

We start by noticing that 
\begin{equation*}
{\tilde S}_n ( \theta)  = \Vert {\tilde \theta} \Vert_e^{-1} \langle {\tilde \theta} , X
\rangle_{{\mathbb{R}}^n}  =  \Vert {\tilde \theta} \Vert_e^{-1} \langle {\tilde \theta} , {\tilde X}
\rangle_{{\mathbb{R}}^n} 
\end{equation*}
where ${\tilde X} = ( X_1 - {\bar X}_n, \cdots, X_n - {\bar X}_n)^t$ and ${%
\tilde \theta} = ( \theta_1 - {\bar \theta}_n, \cdots, \theta_n - {\bar \theta}_n)^t$.

Note that $\sum_{i=1}^n {\tilde \theta}_i = \sum_{i=1}^n {\tilde X}_i = 0$. Let 
$(u_i)_{1 \leq i \leq n}$ be the vectors of ${\mathbb{R}}^n$  defined as
follows: 
\begin{equation*}
u_1= \Big (\frac{1}{\sqrt{2}},-\frac{1}{\sqrt{2}}, 0, \cdots, 0 \Big )^t \,
, \, u_n= \Big (\frac{1}{\sqrt{n}},\frac{1}{\sqrt{n}}, \cdots, \frac{1}{%
\sqrt{n}} \Big )^t
\end{equation*}
and, for $2 \leq k \leq n-1$, 
\begin{equation*}
u_k = \Big (\underbrace{\frac{1}{\sqrt{k(k+1)}}, \cdots, \frac{1}{\sqrt{%
k(k+1)}}}_{k \, times} , - \frac{k}{\sqrt{k(k+1)}}, \underbrace{0, \cdots, 0}%
_{n-k-1 \, times}\Big )^t \, .
\end{equation*}
Note that $(u_i)_{1 \leq i \leq n}$ is an orthonormal basis of ${\mathbb{R}}%
^n$ and that $(u_1, u_2, \cdots, u_n)$ is the change-of-basis matrix from the basis $(u_i)_{1 \leq i
\leq n}$ to the canonical basis $%
(e_i)_{1 \leq i \leq n}$ of ${\mathbb{R}}^n$.  Since $(e_i)_{1 \leq i \leq
n} $ and $(u_i)_{1 \leq i \leq n}$ are both orthonormal bases of ${\mathbb{R}%
}^n $, $(u_1, u_2, \cdots, u_n)$ is an orthonormal matrix and the change-of-basis matrix $B$ from $%
(e_i)_{1 \leq i \leq n}$ to $(u_i)_{1 \leq i \leq n}$ satisfies $B=(u_1, u_2, \cdots, u_n)^t$. Hence ${ \theta} = \sum_{i=1}^{n} (B \theta)_i u_i$ and ${X} =
\sum_{i=1}^{n} (B X)_i u_i$. 
Since ${\tilde \theta}$ (resp. ${\tilde X}$) is the orthogonal  projection of ${ \theta} $ (resp. $X$) on the space generated by $(u_1, \ldots, u_{n-1})$, 
${\tilde \theta} = \sum_{i=1}^{n-1} (B \theta)_i u_i$ and ${\tilde X} =
\sum_{i=1}^{n-1} (B X)_i u_i$.  It follows that 
\begin{equation*}
\langle {\tilde \theta} , {\tilde X} \rangle_{{\mathbb{R}}^n} =
\sum_{i=1}^{n-1} (B \theta )_i (B X )_i \, .
\end{equation*}
In addition, 
\begin{equation}  \label{equalitynorms}
\Vert {\tilde \theta} \Vert_e^2 = \sum_{i=1}^{n-1} (B \theta )^2_i  \ \text{ and } \ 
\Vert {\tilde X} \Vert_e^2 = \sum_{i=1}^{n-1} (B X )^2_i \, .
\end{equation}
So, overall, 
\begin{equation*}
{\tilde S}_n ( \theta) =  \Vert {\tilde \theta} \Vert_e^{-1} \langle {\tilde \theta} , {\tilde X}
\rangle_{{\mathbb{R}}^n}  =:  \langle {\hat \theta} , Y \rangle_{{\mathbb{R}}^{n-1}}
\end{equation*}
where ${\hat \theta}= ( (B \theta )_1, \cdots, (B \theta )_{n-1})^t / \Vert {\tilde \theta}
\Vert_e$ and $Y= ( (BX )_1, \cdots, ( BX )_{n-1})^t$. Note that since $\theta$ has the same law as $\Vert { \xi} \Vert_e^{-1} \xi $ where $\xi$ is a
standard Gaussian vector (i.e. with ${\mathcal{N}} (0, I_n)$ distribution) and $B$ is an
orthonormal matrix, ${\hat \theta}$ is uniformly distributed on the sphere $S^{n-2}$.

Next, setting 
 \begin{equation*}
X^*_{k} = \sqrt{\frac{k}{k+1} } \big ( {\bar X}_k - X_{k+1} \big ) \, ,
\end{equation*}
note that 
\begin{equation*} \label{firstwriting}
\langle  {\hat \theta} , Y \rangle_{{\mathbb{R}}^{n-1}} = \sum_{k=1}^{n-1} {\hat \theta}_k X^*_{k}
\, .
\end{equation*}
Hence by interchanging the sums, it follows that 
\begin{equation*}
{\tilde S}_n ( \theta)= \langle  {\hat \theta} , Y \rangle_{{\mathbb{R}}^{n-1}} = \sum_{\ell=1}^{n}
X_\ell {\ \theta}^*_{\ell} \, .
\end{equation*}
where 
\begin{equation}  \label{deftildetheta}
{ \theta}^*_{n} = - \frac{\sqrt{n-1}}{\sqrt{n}} {\hat \theta}_{n-1} \text{ and } {%
\ \theta}^*_{\ell} = {\ \tilde \theta}_{\ell}- \frac{\sqrt{\ell-1}}{\sqrt{%
\ell}} {\hat \theta}_{\ell-1} \, , 1 \leq \ell \leq n-1 \, ,
\end{equation}
with ${ \tilde \theta}_{\ell} = \sum_{v=\ell}^{n-1} {\hat \theta}_v/ \sqrt{v(v+1)}$
and the convention that ${\hat \theta}_{0} =0$.

\medskip

To summarize ${\tilde S}_n ( \theta) $ can be viewed either as the projection of $(X_1^*,
\ldots, X^*_{k} )$ on the sphere $S_{n-2}$ (however $(X_k^*)_{k \in {\mathbb{%
Z}}}$ is not anymore a sequence of martingale differences) or as the
weighted sums $\sum_{\ell=1}^{n} X_\ell {\ \theta}^*_{\ell} $. Even if $%
(\theta_1^*, \ldots, \theta^*_{n} )$ is not uniformly distributed on the
sphere, we shall use both ways of writings ${\tilde S}_n ( \theta) $ to prove the upper bound \eqref{concluBEmartthm2}.  The proof uses similar arguments as those developed to prove the upper bound \eqref{concluBEmartthm} with substantial modifications that we describe below. 

\medskip

Using the fact that the $X_{k}$'s are uncorrelated, we have that ${\mathbb{E}} ( X^*_{k} )^2 = 1$ for any $k \geq 1$. Then, recalling that 
${\hat \theta}$ is uniformly distributed on the sphere $S^{n-2}$, Lemma 5.2 in \cite{BCG18}  applied with $p=2$ gives: 
for all $t \in {\mathbb R}$, 
\begin{equation}  \label{firstappliLma52}
{\mathbb{E}}\Big |{\mathbb{E}}_{| \theta}\big (\mathrm{e}^{\mathrm{i}t \sum_{k=1}^{n-1} {\hat \theta}_k X^*_{k} } \big ) \Big | \ll \frac{m^2_{4,n} ( X^*) + \sigma^2_{4,n} ( X^*) }{n} + \mathrm{e}^{-t^{2}/16} \, ,
\end{equation}
where
\[
\sigma_{4,n} ( X^*) :=  \frac{1}{ \sqrt{n}}\Big \Vert \sum_{k=1}^{n-1} \big ((X^*_{k} )^2 - {\mathbb{E}} (X^*_{k} )^2 %
\big ) \Big \Vert_2 \  \text{ and } \  m_{4,n} ( X^*)  := \frac{1 }{ \sqrt{n} } \Big ( \E \Big | \sum_{k=1}^{n-1}  X^*_{k}  Y^*_{k}  \Big |^4 \Big )^{1/4} \, , 
\]
with $Y^*=(Y^*_{1} , \cdots, Y^*_{n-1} )$  an independent copy of $X^*= (X^*_{1} , \cdots, X^*_{n-1} )$. But, by independence between $X^*$ and $Y^*$, 
 \[
m_{4,n} (X^*) \leq  \Big (   \sup_{\hat \theta \in S_{n-2}} {\mathbb{E%
}}_{| \hat \theta} \Big ( \sum_{k=1}^{n-1}  \hat \theta_k X^*_{k} \Big )^4 \Big )^{1/4}  \frac{1 }{ \sqrt{n} } \Big ( \E \big | \sum_{k=1}^{n-1}  (Y^*_{k} )^2 \big |^2 \Big )^{1/4}  \, .
\]
But
\[
\Big \Vert \sum_{k=1}^{n-1}  (Y^*_{k} )^2   \Big \Vert_2 =  \Big  \Vert  \sum_{k=1}^{n-1}  (X^*_{k} )^2  \Big  \Vert_2  \leq   \sum_{k=1}^{n-1} \Vert  X^*_{k}  \Vert_4^2 \, .
\]
Using that $\Vert  X^*_{k}  \Vert_4 \leq   2 \Vert  X_0 \Vert_4 $, we then get
 \begin{equation} \label{boundm4n}
m_{4,n} (X^*) \leq 2 \Vert  X_0 \Vert_4   \Big (   \sup_{\hat \theta \in S_{n-2}} {\mathbb{E%
}}_{| \hat  \theta} \Big ( \sum_{k=1}^{n-1} \hat  \theta_k X^*_{k} \Big )^4 \Big )^{1/4}  \, .
\end{equation}
Next, we shall prove that
 \begin{equation}  \label{firstaimLma52}
\Big \Vert \sum_{k=1}^{n-1} \big ((X^*_{k} )^2 - {\mathbb{E}} (X^*_{k} )^2 %
\big ) \Big \Vert_2^2 \ll n \  \text{ and } \  \sup_{ \hat \theta \in S_{n-2}} {\mathbb{E%
}}_{| \hat  \theta} \Big ( \sum_{k=1}^{n-1}  \hat \theta_k X^*_{k} \Big )^4 \ll ( \log n)^2 \, . 
\end{equation}
Starting from \eqref{firstappliLma52} and taking into account the upper bounds  \eqref{boundm4n} and \eqref{firstaimLma52}, it will follow that  for all $T\geq T_{0} \geq 1$ and all $n \geq 2$, 
\begin{equation*}
\int_{T_{0}}^{T}t^{-1}{\mathbb{E}}\Big | {\mathbb{E}}_{| \theta} \big (\mathrm{e}^{\mathrm{i} t\sum_{k=1}^{n-1} {\hat \theta}_k X^*_{k}  } \big )\Big | dt\ll  \log (T/T_{0})\frac { \log n
}{n} +\mathrm{e}^{-T_{0}^{2}/16}\, . \label{B1forregression}
\end{equation*}%
Selecting $T_0= 4 \sqrt{ \log n}$ and $T = n$, we will derive 
\begin{equation*}  \label{B3regression0}
{\mathbb{E}}\big (\kappa _{\theta }(P_{ {\tilde S}_{n}(\theta )},P_{G})\big ) \ll \int_0^{T_0} {\mathbb{E}} %
\Big | {\mathbb{E}}_{| \theta}\Big ( \mathrm{e }^{\mathrm{i}
\sum_{k=1}^{n-1} {\hat \theta}_k X^*_{k} } \Big ) - \mathrm{e}^{-t^2/2} \Big |%
\frac{ dt }{t} + \frac{(\log n)^2}{n} \, ,
\end{equation*}
in place of \eqref{B2} of Proposition \ref{smoothingLemma}. 
Therefore, using that $\sum_{k=1}^{n-1} {\hat \theta}_k X^*_{k} = \sum_{k=1}^{n} {%
\ \theta}^*_k X_k$, where the $\theta_k^*$'s are defined in %
\eqref{deftildetheta}, it will suffice to prove that 
\begin{equation}  \label{B3regression}
\int_0^{T_0} {\mathbb{E}} \Big | {\mathbb{E}}_{| \theta }\Big ( \mathrm{e }^{%
\mathrm{i} \sum_{k=1}^{n} {\theta}^*_k X_{k} } \Big ) - \mathrm{e}^{-t^2/2} %
\Big |\frac{ dt }{t} \ll \frac{(\log n)^{2} }{n} \, , 
\end{equation}
to get the upper bound \eqref{concluBEmartthm2}. 

Let us start by proving \eqref{firstaimLma52}. Note first that 
\begin{multline*}
\Big \Vert \sum_{k=1}^{n-1} \big ((X^*_{k} )^2 - {\mathbb{E}} (X^*_{k} )^2 %
\big ) \Big \Vert_2 \leq \Big \Vert \sum_{k=1}^{n-1} \frac{k}{k+1} \big ( %
X^2_{k+1} - {\mathbb{E}} ( X^2_{k+1}  ) \big ) \Big \Vert_2 \\
+ 2 \Big \Vert \sum_{k=1}^{n-1} \frac{k}{k+1} {\bar X}_k X_{k+1} \Big \Vert%
_2 +  2 \sum_{k=1}^{n-1} \Vert {\bar X}_k \Vert^2_4 \, .
\end{multline*}
Since $(X_k)_{k \geq 0}$ is a sequence of martingale differences in ${%
\mathbb{L}}^4$, applying \eqref{BI} with $p=4$, we get 
\begin{equation} \label{appliBI}
\sum_{k=1}^{n-1} \Vert {\bar X}_k \Vert^2_4 \ll \Vert X_0 \Vert_4^2
\sum_{k=1}^{n-1} k^{-1} \ll \Vert X_0 \Vert_4^2 (\log n)
\end{equation}
and 
\begin{equation*}
\Big \Vert \sum_{k=1}^{n-1} \frac{k}{k+1} {\bar X}_k X_{k+1} \Big \Vert_2^2
= \sum_{k=1}^{n-1} \frac{k^2}{(k+1)^2} \Vert {\bar X}_k X_{k+1} \Vert_2^2
\leq \Vert X_0 \Vert_4^2 \sum_{k=1}^{n-1} \Vert {\bar X}_k \Vert^2_4 \ll
\Vert X_0 \Vert_4^4 (\log n) \, .
\end{equation*}
On another hand, 
\begin{equation*}
\Big \Vert \sum_{k=1}^{n-1} \frac{k}{k+1} \big ( X^2_{k+1} - {\mathbb{E}} (
X^2_{k+1}) \big ) \Big \Vert^2_2 \leq 2 n \sum_{v=0}^{n-1} \gamma_{2,2} (v)
\, .
\end{equation*}
All the above considerations prove that the first part of condition %
\eqref{firstaimLma52} is satisfied as soon as $\sum_{v \geq 0} \gamma (v) <
\infty$. To prove its second part, we first write that 
\begin{equation*}
{\mathbb{E}}_{| \hat  \theta} \Big ( \sum_{k=1}^{n-1} {\hat \theta_k} X^*_{k} \Big )^4
\leq 2^3 {\mathbb{E}}_{| \hat \theta} \Big ( \sum_{k=1}^{n-1} \frac{\sqrt{k}}{%
\sqrt{k+1}}  {\hat \theta_k}  X_{k+1} \Big )^4 + 2^3 {\mathbb{E}}_{| \hat \theta} \Big ( %
\sum_{k=1}^{n-1} \frac{\sqrt{k}}{\sqrt{k+1}}  {\hat \theta_k}  {\bar X}_k \Big )^4 \,
.
\end{equation*}
Burkholder's inequality \eqref{BI} implies that 
\begin{equation*}
{\mathbb{E}}_{| \hat \theta} \Big ( \sum_{k=1}^{n-1} \frac{\sqrt{k}}{\sqrt{k+1}}
 {\hat \theta_k}  X_{k+1} \Big )^4 \ll \Vert X_0 \Vert_4^4 \Big ( \sum_{k=1}^{n-1} 
\frac{k}{k+1}  {\hat \theta_k}^2 \Big )^2 \ll \Vert X_0 \Vert_4^4 \, .
\end{equation*}
Next,
\begin{equation*}
{\mathbb{E}}_{|  \hat \theta} \Big ( \sum_{k=1}^{n-1} \frac{\sqrt{k}}{\sqrt{k+1}}
 {\hat \theta_k}  {\bar X}_k \Big )^4 \leq  \Big ( \sum_{k=1}^{n-1}  \frac{k}{k+1} {\hat \theta_k} ^2 \Big )^2 {\mathbb{E}} \Big ( \sum_{k=1}^{n-1} {\bar
X}_k ^2 \Big )^2 \leq \Big ( \sum_{k=1}^{n-1} \Vert {\bar X}_k \Vert_4^2 %
\Big )^2 \, .
\end{equation*}
Hence, by \eqref{appliBI}, 
\[
{\mathbb{E}}_{|\hat  \theta} \Big ( \sum_{k=1}^{n-1} \frac{\sqrt{k}}{\sqrt{k+1}}
 {\hat \theta_k}  {\bar X}_k \Big )^4\ll \Vert {\ X}_0 \Vert_4^4 (\log n)^2 \, .
 \]
So, overall, the second part of \eqref{firstaimLma52} is proved.

\smallskip

It remains to show that \eqref{B3regression} is satisfied. With this aim, we
consider $(Y_i( {\theta}^*) )_{1 \leq i \leq n}$ a sequence of random
variables that are independent for any fixed $\theta$, independent of $(X_i
)_{i \geq 1}$, and such that, for each $i \geq 1$, the conditional law of $%
Y_i ( {\theta}^*)$ given $\theta$ takes 2 values and is such that ${\mathbb{E%
}}_{| \theta} ( Y_i( {\theta}^*) ) =0$, ${\mathbb{E}}_{| \theta} ( Y^2_i( {\theta}^*) )
= ( {\theta}^*_i)^2 {\mathbb{E}} ( X_0^2)$, ${\mathbb{E}}_{| \theta} (
Y^3_i(  {\theta}^*) ) =\beta_{i,3} ( {\theta}^*)$ and ${\mathbb{E}}_{|
\theta} ( Y^4_i( {\theta}^*) ) =\beta_{i,4} ( {\theta}^*)$ where $%
\beta_{i,3} ( {\theta}^*)$ (resp. $\beta_{i,4} ( {\theta}^*)$) is defined by %
\eqref{defbeta3k} (resp. \eqref{defbeta4k}) with ${\theta}^*$ replacing ${%
\theta}$. Recall that this is always possible according to Fact \ref%
{defYtheta} (when $ {\theta}^*_i=0$, we set $\beta_{i,4} ( {\theta}^*) =0$ and $Y_i ( {\theta}^*)=0$).

To show that \eqref{B3regression} is satisfied, we shall prove that for all $n \geq 2$, 
\begin{equation}  \label{B3regression1}
\int_0^{T_0} {\mathbb{E}} \Big ( \Big | {\mathbb{E}}_{| \theta} \Big ( 
\mathrm{e }^{\mathrm{i} \sum_{k=1}^{n} {\theta}^*_k X_k } \Big ) - {\mathbb{E%
}}_{| \theta} \Big ( \mathrm{e }^{\mathrm{i} \sum_{k=1}^{n} Y_k ( {\theta}%
^*) } \Big ) \Big | \Big )\frac{ dt }{t} \ll \frac{(\log n)^{2} }{n} \, ,
\end{equation}
\begin{equation}  \label{B3regression2}
\int_0^{T_0} {\mathbb{E}} \Big ( \Big | {\mathbb{E}}_{| \theta} \Big ( 
\mathrm{e }^{\mathrm{i} \sum_{k=1}^{n} Y_k ( {\theta}^*) } \Big ) - \mathrm{e%
}^{-t^2\sum_{k=1}^{n} ( {\theta}^*_{k})^2/2} \Big | \Big )\frac{ dt }{t} \ll 
\frac{(\log n)^{2} }{n} \, ,
\end{equation}
and 
\begin{equation}  \label{B3regression3}
\int_0^{T_0} {\mathbb{E}} \big | \mathrm{e}^{-t^2\sum_{k=1}^{n} ({\theta}%
^*_{k})^2/2} - \mathrm{e}^{-t^2/2} \big | \frac{ dt }{t} \ll\frac{(\log
n)^{2} }{n} \, .
\end{equation}
To prove \eqref{B3regression1} and \eqref{B3regression2}, we proceed as in
the proof of Theorem \ref{BEmartthm}. Then starting from \eqref{concludec}
with ${\theta}^*$ instead of $\theta$, the estimate \eqref{B3regression1}
will follow if one can prove that, for any $m \geq 1$, there exists a
positive constant $c_m$ such that, any $1 \leq k \leq n$, 
\begin{equation}  \label{PR1}
{\mathbb{E}} (| {\theta}^*_k|^m) \leq c_m n^{-m/2} \, ,
\end{equation}
and 
\begin{equation}  \label{PR2}
\sum_{k=1}^{n} {\mathbb{E}} ( \beta_{k,4} ( {\theta}^* ) ) \ll n^{-1} \, .
\end{equation}
On another hand, analyzing the proof of \eqref{B3indepart}, we infer that %
\eqref{B3regression2} holds provided \eqref{PR2} is satisfied and 
\begin{equation}  \label{PR2*}
{\mathbb{E}} \Big | \sum_{k=1}^{n} \beta_{k,3} ({\theta}^* ) \Big | \ll
n^{-1} \, , \, \sum_{k=1}^{n} {\mathbb{E}} ( \beta^2_{k,3} ({\theta}^* ) )
\ll n^{-1} \, \text{ and } \, \max_{1 \leq k \leq n} {\mathbb{E}} (
\beta^{m}_{k,4} ({\theta}^* ) ) \leq c_m n^{-2m} \, ,
\end{equation}
for any $m \geq 1$ (above $c_m$ is a constant not depending on $n$). 
Finally, to prove \eqref{B3regression3}, we note that 
\begin{equation*}
{\mathbb{E}} \Big | \mathrm{e}^{-t^2 \sum_{k=1}^{n} ({\theta}^*_{k})^2 /2} %
 - \mathrm{e}^{-t^2 /2} \Big | \leq \frac{t^2}{2} {\mathbb{E}} \Big | %
\sum_{k=1}^{n} ({\theta}^*_k)^2 -1 \Big | \, .
\end{equation*}
Hence to prove \eqref{B3regression3}, it suffices to show that 
\begin{equation}  \label{PR3}
{\mathbb{E}} \Big | \sum_{k=1}^{n} ({\theta}^*_{k})^2 - 1 \Big | \ll n^{-1}
\log n \, .
\end{equation}

We start by proving \eqref{PR1}. For any positive integer $k \leq n$ and any $m \geq 1$, note first that 
\begin{equation*}
{\mathbb{E}} ( | {\theta}^*_k|^m) \leq 2^{m-1} {\mathbb{E}} ( | {  \hat \theta}%
_{k-1}|^m) + 2^{m-1} {\mathbb{E}} \Big ( \Big | \sum_{\ell=k}^{n-1} \frac{
{  \hat \theta}_\ell }{\sqrt{\ell ( \ell+1)}}\Big |^m \Big )
\end{equation*}
Since ${  \hat \theta}$ is uniformly distributed on the sphere $S^{n-2}$, ${\mathbb{E}%
} ( | {  \hat \theta}_{k-1}|^m) \ll c_m n^{-m/2}$. In addition $({  \hat \theta}_1,
\ldots, {  \hat \theta}_{n-1}) =^{\mathcal{D}} ( \xi_1, \ldots, \xi_{n-1}) / \Vert
\xi \Vert_e $ where $(\xi_i)_{i \geq 1}$ is a sequence of iid ${\mathcal{N%
}} (0,1)$-distributed r.v.'s and $\Vert \xi \Vert_e^2 = \sum_{i=1}^{n-1} \xi_i^2$. Hence 
\begin{multline*}
{\mathbb{E}} \Big ( \Big | \sum_{\ell=k}^{n-1} \frac{ {  \hat \theta}_\ell }{\sqrt{%
\ell ( \ell+1)}}\Big |^m \Big ) \ll n^{-m/2} {\mathbb{E}} \Big ( \Big | %
\sum_{\ell=k}^{n-1} \frac{ \xi_\ell }{\sqrt{\ell ( \ell+1)}}\Big |^m \Big )
\\
+ {\mathbb{E}}^{1/2} \Big ( \Big | \sum_{\ell=k}^{n-1} \frac{ {  \hat \theta}_\ell }{ 
\sqrt{\ell ( \ell+1)}}\Big |^{2m} \Big ) \Big ( {\mathbb{P}} \big ( \Vert \xi
\Vert_e^2 < n/2 \big ) \Big )^{1/2} \, .
\end{multline*}
By the Burkholder inequality \eqref{BI}, for any $m \geq 1$, 
\begin{multline*}
 \Big \Vert  \sum_{\ell=k}^{n-1} \frac{ \xi_\ell }{\sqrt{\ell
( \ell+1)}}\  \Big \Vert_m \leq  \Big \Vert  \sum_{\ell=k}^{n-1} \frac{ \xi_\ell }{\sqrt{\ell
( \ell+1)}}\  \Big \Vert_{m \vee 2} \\ \leq  \sqrt{m \vee 2 -1 }  \Vert  \xi_0 \Vert_{m \vee 2} \Big | \sum_{\ell=k}^{n-1} \frac{1 }{\ell (
\ell+1)}\Big |^{1/2}  \leq c_m k^{-1/2} \, .
\end{multline*}
On another hand, using that ${\mathbb{P}} \big ( \Vert \xi \Vert_e^2 < n/2 %
\big )  \leq e^{-cn}$, for some $c>0$, and the fact that ${  \hat \theta}_k \leq 1$,
we get 
\begin{equation*}
{\mathbb{E}}^{1/2} \Big ( \Big | \sum_{\ell=k}^{n-1} \frac{ {  \hat \theta}_\ell }{ 
\sqrt{\ell ( \ell+1)}}\Big |^{2m} \Big ) \Big ( {\mathbb{P}} \big ( \Vert \xi
\Vert_e^2 < n/2 \big ) \Big )^{1/2} \ll n^{-m/2} k^{-m/2} \, .
\end{equation*}
So, overall, 
\begin{equation}  \label{interthetatildem}
{\mathbb{E}} (| {  \tilde \theta}_k |^m )= {\mathbb{E}} \Big ( \Big | \sum_{\ell=k}^{n-1} \frac{ {  \hat \theta}_\ell }{\sqrt{%
\ell ( \ell+1)}}\Big |^m \Big ) \ll n^{-m/2} k^{-m/2}\, ,
\end{equation}
ending the proof of  \eqref{PR1}. 

We turn now to the proof of the last part of %
\eqref{PR2*} that will also imply \eqref{PR2}. Setting $a(k) ={\mathbb{E}}
( X_{0} X_k^2)$, note first that 
\begin{equation*}
\beta_{k,4} ( {\theta}^*) \ll {  \hat \theta}_{k-1}^4 + \Big ( \sum_{\ell=k}^{n-1} 
\frac{ {  \hat \theta}_\ell }{\sqrt{\ell ( \ell+1)}} \Big )^4 + \Big ( \sum_{\ell
=1}^{k-1} \frac{ \sqrt{\ell-1}}{ \sqrt{\ell}} {  \hat \theta}_{\ell-1} a(k-\ell) \Big
)^4 + \Big ( \sum_{\ell =1}^{k-1} {\tilde \theta}_{\ell} a(k - \ell) \Big )%
^4 \, .
\end{equation*}
Proceeding as in the proof of \eqref{interthetatildem},  since by assumption $%
\sum_{k \geq 1} a^2(k) < \infty$, we infer that
\begin{equation*}
\max_{1 \leq k \leq n}{\mathbb{E}} \Big ( \sum_{\ell =1}^{k-1} \frac{ \sqrt{%
\ell-1}}{ \sqrt{\ell}} {  \hat \theta}_{\ell-1} a(k-\ell) \Big )^{4m} \ll n^{-2m} \, ,
\end{equation*}
which combined with \eqref{interthetatildem} and the fact that $\max_{1 \leq
k \leq n}{\mathbb{E}} ( {  \hat \theta}_{k-1}^{4m}) \ll n^{-2m}$, implies that, for
any $m \geq 1$, 
\begin{equation*}
\max_{1 \leq k \leq n} {\mathbb{E}} (\beta^m_{k,4} ( { \theta}^*) ) \ll
n^{-2m} + \max_{1 \leq k \leq n} {\mathbb{E}} \Big ( \sum_{\ell =1}^{k-1} {%
\tilde \theta}_{\ell} a(k - \ell) \Big )^{4m} \, .
\end{equation*}
Next, proceeding again as in the proof of \eqref{interthetatildem}, for any $k \leq
n$, we get 
\begin{multline}  \label{PR2-3}
{\mathbb{E}} \Big ( \sum_{\ell =1}^{k-1} {\tilde \theta}_{\ell} a(k - \ell) %
\Big )^{4m} = {\mathbb{E}} \Big ( \sum_{v =1}^{n-1} \frac{ { \hat \theta}_{v}}{%
\sqrt{v (v+1)}} \sum_{\ell =1}^{(k-1) \wedge v } a(k - \ell) \Big )^{4m} \\
\ll n^{-2m} + n^{-2m} \Big [ \sum_{v =1}^{n-1} \frac{ 1}{v (v+1)} \Big ( %
\sum_{\ell =1}^{(k-1) \wedge v } a(k - \ell) \Big )^2 \Big ]^{2m} \ll
n^{-2m} \Big ( 1 + \sum_{\ell = 1}^{n-1} \gamma (v) \Big )^{4m} \, .
\end{multline}
So, overall, for any $m \geq 1$, 
\begin{equation*}
\max_{1 \leq k \leq n}{\mathbb{E}} (\beta^m_{k,4} ( { \theta}^*) ) \ll
n^{-2m} \Big ( 1 + \sum_{\ell = 1}^{n-1} \gamma (v) \Big )^{4m} \, ,
\end{equation*}%
which proves the last part of \eqref{PR2*} (and then \eqref{PR2}) since $%
\sum_{v \geq 1} \gamma (v) < \infty$. 

We prove now the second part of  \eqref{PR2*}. By Young's
inequality, 
\begin{equation*}
\sum_{k=1 }^{n} {\mathbb{E}} \big ( \beta^2_{k,3} ( { \theta}^*) \big ) \ll
\sum_{k=1 }^{n-1} {\mathbb{E}} ( { { \theta}_k^*}^6 ) + \sum_{k=1 }^{n} {%
\mathbb{E}} \Big ( \sum_{\ell =1}^{k-1} \frac{ \sqrt{\ell-1}}{ \sqrt{\ell}}
a(k - \ell) {  \hat \theta}_{\ell-1} \Big )^6 + \sum_{k=1 }^{n} {\mathbb{E}} \Big ( %
\sum_{\ell =1}^{k-1} {\tilde \theta}_{\ell} a(k - \ell) \Big )^6 \, .
\end{equation*}
By \eqref{PR1}, the first term in the right-hand side is bounded by a
constant times $n^{-2}$. In addition, by the same arguments used to prove %
\eqref{interthetatildem} and \eqref{PR2-3} (with $m=3/2$), we infer that 
\begin{equation}  \label{PR2bis}
\sum_{k=1 }^{n} {\mathbb{E}} \Big ( \sum_{\ell =1}^{k-1} \frac{ \sqrt{\ell-1}%
}{ \sqrt{\ell}} a(k - \ell) {  \hat \theta}_{\ell-1} \Big )^6 \ll n^{-3} \sum_{k=1
}^{n} \Big ( \sum_{\ell =1}^{k-1} a^2(k-\ell) \Big )^3 + n^{-2} \ll n^{-2}
\, ,
\end{equation}
and 
\begin{equation*}
\sum_{k=1 }^{n} {\mathbb{E}} \Big ( \sum_{\ell =1}^{k-1} {\tilde \theta}%
_{\ell} a(k - \ell) \Big )^6 \ll n^{-3} \Big ( 1 +
\sum_{\ell = 1}^{n-1} \gamma (v) \Big )^6 \ll n^{-2} \, .
\end{equation*}
So, overall, the second part of \eqref{PR2*} is satisfied. 

We turn to the
proof of \eqref{PR3}. With this aim, we note that $\sum_{k=1}^{n-1} {\hat \theta%
}_{k}^2 =1$. Hence 
\begin{equation*}
\sum_{k=1}^{n} { { \theta}_{k}^*}^2 - 1 = \sum_{\ell=1}^{n-1} ({ { \tilde
\theta}_{\ell}})^2 - 2 \sum_{\ell=2}^{n-1} \frac{\sqrt{\ell-1}}{\sqrt{\ell}} {%
\tilde \theta}_{\ell} {  \hat \theta}_{\ell-1} - \sum_{\ell=1}^{n-1} \frac{1}{\ell +1}
{  \hat \theta}^2_{\ell} \, .
\end{equation*}
Taking the expectation of the absolute values of the above quantity, and
considering  \eqref{interthetatildem}, we then infer that %
\eqref{PR3} holds. Indeed to deal with the second term in the right-hand side,  we can use the same arguments used to prove %
\eqref{interthetatildem} and the fact that 
\begin{multline*}
{\mathbb{E}} \Big | \sum_{\ell=2}^{n-1} \frac{\sqrt{\ell-1}}{\sqrt{\ell}}
\xi_{\ell-1} \sum_{v=\ell}^{n-1} \frac{\xi_v}{ \sqrt{v(v+1)}} \Big |^2 = {%
\mathbb{E}} \Big | \sum_{v=2}^{n-1} \frac{\xi_v}{ \sqrt{v(v+1)}}
\sum_{\ell=2}^{v} \frac{\sqrt{\ell-1}}{\sqrt{\ell}} \xi_{\ell-1} \Big |^2 \\
= \sum_{v=2}^{n-1} \frac{1}{ v(v+1)} {\mathbb{E}} \Big ( \sum_{\ell=2}^{v} 
\frac{\sqrt{\ell-1}}{\sqrt{\ell}} \xi_{\ell-1} \Big )^2 \leq
\sum_{v=2}^{n-1} \frac{1}{ v+1} \leq \log n \, .
\end{multline*}
to derive that 
\begin{equation}  \label{PR3-P1}
{\mathbb{E}} \Big | \sum_{\ell=2}^{n-1} \frac{\sqrt{\ell-1}}{\sqrt{\ell}}
{  \hat \theta}_{\ell-1} \sum_{v=\ell}^{n-1} \frac{{  \hat \theta}_v}{ \sqrt{v(v+1)}} \Big |^2
\ll n^{-2} \log n \, .
\end{equation}

To end the proof of the theorem, it remains to prove the first part of %
\eqref{PR2*}. Recall that 
\begin{equation*}
\sum_{k=1}^n \beta_{k,3} ( { \theta}^*) = {\mathbb{E}} ( X_0^3)
\sum_{k=1}^n {\theta_{k}^* }^3 + 3\sum_{k=1}^n {\theta_{k}^* }^2 \sum_{\ell
=1}^{k-1} {\theta^*_{\ell}} a(k - \ell) \, .
\end{equation*}
Note that 
\begin{multline*}
\Big | \sum_{k=1}^n {\theta_{k}^*}^3 \Big | \ll \Big | \sum_{k=1}^{n-1} ({\tilde \theta_{k}})^3 \Big | + \Big | \sum_{k=1}^n \Big ( \frac{k-1}{k} \Big )%
^{3/2}{{  \hat \theta}^3_{k-1}} \Big | \\ + \Big | \sum_{k=1}^{n-1} ( { \tilde \theta_{k}})^2 \Big ( \frac{k-1}{k} \Big )^{1/2}{{  \hat \theta}_{k-1}} \Big | + \Big | %
\sum_{k=1}^{n-1} { \tilde \theta_{k}} \frac{k-1}{k} {{  \hat \theta}^2_{k-1}} \Big | \, .
\end{multline*}
By  \eqref{interthetatildem}, 
\begin{equation*}
{\mathbb{E}} \Big | \sum_{k=1}^{n-1} ({ \tilde \theta_{k}})^3 \Big | \leq
\sum_{k=1}^{n-1} {\mathbb{E}} (| ({\tilde \theta_{k}})^3 |) \ll n^{-3/2}
\sum_{k=1}^n k^{-3/2} \ll n^{-3/2}
\end{equation*}
and 
\begin{equation*}
{\mathbb{E}} \Big | \sum_{k=1}^{n-1}  ({ \tilde \theta_{k}})^2 \Big ( \frac{k-1}{k}
\Big )^{1/2}{ {  \hat \theta}_{k-1}} \Big | \leq \sum_{k=1}^{n-1} {\mathbb{E}}^{1/2} ({\tilde \theta_{k}}^4 ) {\mathbb{E}}^{1/2} ({  \hat \theta}^2_{k-1} ) \ll n^{-3/2}
\sum_{k=1}^{n-1} k^{-1} \ll n^{-3/2}\log n \, .
\end{equation*}
Next, using the same arguments used for proving \eqref{interthetatildem}, we infer that 
\begin{equation*}
{\mathbb{E}} \Big | \sum_{k=1}^n \Big ( \frac{k-1}{k} \Big )^{3/2}{\
{  \hat \theta}^3_{k-1}} \Big | \ll n^{-1} + n^{-3/2} \sqrt{ {\mathbb{E}} \Big | %
\sum_{k=1}^n \Big ( \frac{k-1}{k} \Big )^{3/2}{\ \xi^3_{k-1}} \Big |^2 } \ll
n^{-1} \, .
\end{equation*}
On another hand, by Cauchy-Schwarz's inequality, the fact that ${\mathbb{E}%
} ( | {  \hat \theta}_{k-1}|^m) \ll c_m n^{-m/2}$ and    \eqref{interthetatildem}, we get
\begin{multline*}
{\mathbb{E}} \Big ( \sum_{k=1}^{n-1} { \tilde \theta_{k}} \frac{k-1}{k} {{  \hat \theta}^2_{k-1}} \Big )^2  \leq n  \sum_{k=1}^{n-1}  {\mathbb{E}} \big (   { \tilde \theta_{k}}^2  {{  \hat \theta}^4_{k-1}}\big ) \leq n  \sum_{k=1}^{n-1}  {\mathbb{E}}^{1/2} \big (   { \tilde \theta_{k}}^4  \big )  {\mathbb{E}}^{1/2} \big (    {{  \hat \theta}^8_{k-1}}\big )
\\  \ll n^{-2} \sum_{k=1}^{n-1}  k^{-1} \ll 
n^{-2} (\log n) \, .
\end{multline*}
So, overall, 
\begin{equation}  \label{sumthetaketoile3}
{\mathbb{E}} \Big | \sum_{k=1}^n {\theta_{k}^* }^3 \Big | \ll n^{-1} \, .
\end{equation}
We now give an upper bound for ${\mathbb{E}} \Big | \sum_{k=1}^n {%
\theta_{k}^* }^2 \sum_{\ell =1}^{k-1} {\theta^*_{\ell}} a(k - \ell) \Big |$.
By using  \eqref{interthetatildem}, \eqref{PR2-3} and the fact that $\sum_{v \geq 1}
\gamma (v) < \infty$, we first notice that 
\begin{equation}  \label{PR2etoilep1}
{\mathbb{E}} \Big | \sum_{k=1}^n  ({ \tilde \theta_{k} })^2 \sum_{\ell
=1}^{k-1} { \tilde \theta_{\ell}} a(k - \ell) \Big | \leq \sum_{k=1}^n {%
\mathbb{E}}^{1/2} ( {\tilde \theta_{k} }^4) {\mathbb{E}}^{1/2} \Big ( %
\sum_{\ell =1}^{k-1} {\ \tilde \theta_{\ell}} a(k - \ell) ) \Big )^2 \ll
n^{-3/2} \log n \, .
\end{equation}
On another hand, proceeding as to prove \eqref{boundsumbetak3} and since $\sum_{v \geq 1}
\gamma (v) < \infty$, we get that 
\begin{equation}  \label{PR2etoilep2}
{\mathbb{E}} \Big | \sum_{k=1}^n \frac{k-1}{k} {  \hat \theta}_{k-1}^2 \sum_{\ell
=1}^{k-1} {\ \frac{\sqrt{\ell-1}}{\sqrt{\ell}} {  \hat \theta}_{\ell-1}} a(k - \ell) %
\Big | \ll n^{-1} \, .
\end{equation}
Next 
\begin{equation*}
{\mathbb{E}} \Big | \sum_{k=1}^{n-1} ( { \tilde \theta_{k} })^2 \sum_{\ell
=1}^{k-1} {\ \frac{\sqrt{\ell-1}}{\sqrt{\ell}} {  \hat \theta}_{\ell-1}} a(k - \ell) %
\Big | \leq \sum_{k=1}^{n-1} {\mathbb{E}}^{1/2} ( { \tilde \theta_{k} }^4) {%
\mathbb{E}}^{1/2} \Big ( \sum_{\ell =1}^{k-1} {\ \frac{\sqrt{\ell-1}}{\sqrt{%
\ell}} {  \hat \theta}_{\ell-1}} a(k - \ell) \Big )^2 \, .
\end{equation*}
Using  \eqref{interthetatildem} and proceeding as in the proof of \eqref{PR2bis}, we derive that 
\begin{equation}  \label{PR2etoilep3}
{\mathbb{E}} \Big | \sum_{k=1}^{n-1} ({ \tilde \theta_{k} })^2 \sum_{\ell
=1}^{k-1} {\ \frac{\sqrt{\ell-1}}{\sqrt{\ell}} {  \hat \theta}_{\ell-1}} a(k - \ell) %
\Big | \ll n^{-3/2} \log n \, .
\end{equation}
Next 
\begin{multline*}
\Big | \sum_{k=1}^n \frac{k-1}{k} {{  \hat \theta}^2_{k-1}} \sum_{\ell=1}^{k-1} {%
\tilde \theta}_{\ell} a(k - \ell ) \Big | = \Big | \sum_{k=1}^n \frac{k-1}{k}
{{  \hat \theta}^2_{k-1}} \sum_{v=1}^{n-1} \frac{{  \hat \theta}_{v}}{ \sqrt{v(v+1)}}%
\sum_{\ell=1}^{(k-1) \wedge v}a(k - \ell ) \Big | \\
\leq \Big | \sum_{k=1}^n \frac{k-1}{k} {{  \hat \theta}^2_{k-1}} \sum_{v=1}^{k-1} 
\frac{{  \hat \theta}_{v}}{ \sqrt{v(v+1)}}\sum_{\ell=1}^{ v}a(k - \ell ) \Big | + %
\Big | \sum_{k=1}^{n-1} \frac{k-1}{k} {{  \hat \theta}^2_{k-1}} \sum_{v=k}^{n-1} 
\frac{{  \hat \theta}_{v}}{ \sqrt{v(v+1)}}\sum_{\ell=1}^{k-1}a(k - \ell ) \Big | \, .
\end{multline*}
Using the same arguments used to prove \eqref{interthetatildem}, we infer that 
\begin{multline*}
{\mathbb{E}} \Big | \sum_{k=1}^n \frac{k-1}{k} {{  \hat \theta}^2_{k-1}}
\sum_{v=[k/2] +1}^{k-1} \frac{{  \hat \theta}_{v}}{ \sqrt{v(v+1)}}\sum_{\ell=1}^{
v}a(k - \ell ) \Big | \\
\ll n^{-1} + n^{-3/2} {\mathbb{E}} \Big | \sum_{k=1}^n \frac{k-1}{k} {\
\xi^2_{k-1}} \sum_{v=[k/2] +1}^{k-1} \frac{\xi_{v}}{ \sqrt{v(v+1)}}%
\sum_{\ell=1}^{ v}a(k - \ell ) \Big | \\
\ll n^{-1} + n^{-3/2}  \sum_{\ell=1}^{ n-1} |a(u)|  \sum_{k=2}^n   \frac{ 1}{ \sqrt{k(k-1)}}  + n^{-3/2} \sum_{k=1}^n {\mathbb{E}}^{1/2} \Big ( \sum_{v=[k/2]
+1}^{k-2} \frac{\xi_{v}}{ \sqrt{v(v+1)}}\sum_{\ell=k-v}^{ k-1}a(u) \Big )^2
\\
\ll n^{-1} + n^{-3/2} \sum_{k=1}^n k^{-1} \sqrt{ \sum_{v=[k/2] +1}^{k-2} %
\Big ( \sum_{\ell=k-v}^{ k-1}a(u) \Big )^2 } \ll n^{-1} + n^{-3/2}
\sum_{k=1}^n k^{-1/2} \ll n^{-1} \, .
\end{multline*}
In addition 
\begin{multline*}
{\mathbb{E}} \Big | \sum_{k=1}^n \frac{k-1}{k} {{  \hat \theta}^2_{k-1}}
\sum_{v=1}^{[k/2] } \frac{{  \hat \theta}_{v}}{ \sqrt{v(v+1)}}\sum_{\ell=1}^{ v}a(k -
\ell ) \Big | \\
\ll n^{-3/2} + n^{-3/2} {\mathbb{E}} \Big | \sum_{k=1}^n \frac{k-1}{k} {\
\xi^2_{k-1}} \sum_{v=1}^{[k/2] } \frac{\xi_{v}}{ \sqrt{v(v+1)}}%
\sum_{\ell=1}^{ v}a(k - \ell ) \Big | \\
\ll n^{-3/2} + n^{-3/2} \sum_{k=1}^n k \gamma([k/2]) \ll n^{-3/2} \, .
\end{multline*}
On another hand 
\begin{equation*}
{\mathbb{E}} \Big | \sum_{k=1}^{n-1} \frac{k-1}{k} {{  \hat \theta}^2_{k-1}}
\sum_{v=k}^{n-1} \frac{{  \hat \theta}_{v}}{ \sqrt{v(v+1)}}\sum_{\ell=1}^{k-1}a(k -
\ell ) \Big | = {\mathbb{E}} \Big | \sum_{v=1}^{n-1} \frac{\hat \theta_{v}}{ 
\sqrt{v(v+1)}} \sum_{k=1}^{v} \sum_{\ell=1}^{k-1}a(k - \ell ) \frac{k-1}{k} {%
\hat \theta^2_{k-1}} \Big | \, .
\end{equation*}
By the same arguments as to prove \eqref{interthetatildem}, note that 
\begin{multline*}
{\mathbb{E}} \Big | \sum_{v=1}^{n-1} \frac{{  \hat \theta}_{v}}{ \sqrt{v(v+1)}}
\sum_{k=1}^{v} \sum_{\ell=1}^{k-1}a(k - \ell ) \frac{k-1}{k} {\
{  \hat \theta}^2_{k-1}} \Big | \\
\ll n^{-1} + n^{-3/2} {\mathbb{E}} \Big | \sum_{v=1}^{n-1} \frac{\xi_{v}}{ 
\sqrt{v(v+1)}} \sum_{k=1}^{v} \sum_{\ell=1}^{k-1}a(k - \ell ) \frac{k-1}{k} {%
\ \xi^2_{k-1}} \Big | \, .
\end{multline*}
But, by independence, 
\begin{multline*}
n^{-3/2} {\mathbb{E}}^{1/2} \Big | \sum_{v=1}^{n-1} \frac{\xi_{v}}{ \sqrt{%
v(v+1)}} \sum_{k=1}^{v} \sum_{\ell=1}^{k-1}a(k - \ell ) \frac{k-1}{k} {\
\xi^2_{k-1}} \Big |^2 \\
\leq n^{-3/2} \Big ( \sum_{v=1}^{n-1} \frac{1}{ v(v+1)} {\mathbb{E}} \Big ( %
\sum_{k=1}^{v} \sum_{\ell=1}^{k-1}a(k - \ell ) \frac{k-1}{k} { \xi^2_{k-1}} %
\Big )^2 \Big )^{1/2} \\
\ll n^{-3/2} \Big ( \sum_{v=1}^{n-1} \Big ( \sum_{u=1}^{v-1}\gamma(u) %
\Big )^2 \Big )^{1/2} \ll n^{-1} \, .
\end{multline*}
So, overall, 
\begin{equation}  \label{PR2etoilep4}
{\mathbb{E}}\Big | \sum_{k=1}^n \frac{k-1}{k} {{  \hat \theta}^2_{k-1}}
\sum_{\ell=1}^{k-1} {\tilde \theta}_{\ell} a(k - \ell ) \Big | \ll n^{-1} \,
.
\end{equation}
Moreover, using  \eqref{interthetatildem} and $\sup_{1 \leq \ell \leq n-1}{\mathbb{E}%
} ( | {  \hat \theta}_{\ell}|^m) \ll c_m n^{-m/2}$, we get 
\begin{equation}  \label{PR2etoilep5}
{\mathbb{E}}\Big | \sum_{k=1}^n \frac{\sqrt{k-1}}{\sqrt{k}} {{  \hat \theta}_{k-1}} {\
\tilde \theta_{k}} \sum_{\ell=1}^{k-1} {\tilde \theta}_{\ell} a(k - \ell ) %
\Big | \ll n^{-1} \, ,
\end{equation}
and 
\begin{equation}  \label{PR2etoilep6}
{\mathbb{E}}\Big | \sum_{k=1}^n \frac{\sqrt{k-1}}{\sqrt{k}} {{  \hat \theta}_{k-1}} {\
\tilde \theta_{k}} \sum_{\ell=1}^{k-1} \frac{\sqrt{\ell-1}}{ \sqrt{\ell}}   {  \hat \theta}_{\ell-1} a(k - \ell ) \Big | \ll n^{-1} \, .
\end{equation}
Hence considering the upper bounds \eqref{PR2etoilep1}-\eqref{PR2etoilep6}, we get that 
\begin{equation}  \label{PR2etoileconclu}
{\mathbb{E}} \Big | \sum_{k=1}^n {%
\theta_{k}^* }^2 \sum_{\ell =1}^{k-1} {\theta^*_{\ell}} a(k - \ell) \Big | \ll n^{-1} \, .
\end{equation}
The first part of \eqref{PR2*} then follows by considering the upper bounds %
\eqref{sumthetaketoile3} and \eqref{PR2etoileconclu}. This ends the proof of the upper bound \eqref{concluBEmartthm2}. \qed

\smallskip

\noindent \textbf{Acknowledgements.} This paper was partially supported by the NSF
grant DMS-2054598 and a Taft research support grant.  The authors are  
grateful to the referees for carefully reading our manuscript and for helpful suggestions that significantly improved the presentation of the paper.

\end{document}